\documentclass[final]{siamltex}

\usepackage[usenames]{color}
\usepackage[color]{showkeys}
\definecolor{refkey}{gray}{0.4}
\definecolor{labelkey}{gray}{0.3}

\usepackage[left=1in,top=1in,right=1in,bottom=1in,dvips,letterpaper]{geometry}
\usepackage{setspace,url}
\usepackage[leqno]{amsmath}
\usepackage{amsxtra, amsfonts,amscd,  amssymb, graphicx,subfigure}
\usepackage[algo2e, vlined,ruled]{algorithm2e}
\numberwithin{equation}{section}

\newcommand{\be}{\begin{equation}}
\newcommand{\ee}{\end{equation}}
\newcommand{\ba}{\begin{array}}
\newcommand{\ea}{\end{array}}
\newcommand{\bea}{\begin{eqnarray}}
\newcommand{\eea}{\end{eqnarray}}
\newcommand{\beaa}{\begin{eqnarray*}}
\newcommand{\eeaa}{\end{eqnarray*}}

\newcommand{\half}{\frac{1}{2}}
\newcommand{\br}{\mathbb{R}}

\usepackage[usenames]{color}
\usepackage[normalem]{ulem}

\newcommand{\TV}{\mathrm{TV}}

\newcommand{\etal}{{et al. }}

\newcommand{\oneotwomu}{\frac{1}{2\mu}}
\newcommand{\oneomu}{\frac{1}{\mu}}

\newtheorem{remark}[theorem]{Remark}

\newtheorem{assumption}[theorem]{Assumption}

\begin{document}

\title{Fast Multiple Splitting Algorithms for Convex Optimization}
\date{}
\author{Donald Goldfarb\footnotemark[2] \and Shiqian Ma\footnotemark[2]}

\renewcommand{\thefootnote}{\fnsymbol{footnote}}
\footnotetext[2]{Department of Industrial Engineering and Operations
Research, Columbia University, New York, NY 10027, USA. \quad Email:
\{goldfarb, sm2756\}@columbia.edu. Research supported in
part by NSF Grants DMS 06-06712 and DMS 10-16571, ONR Grant
N00014-08-1-1118 and DOE Grant DE-FG02-08ER25856.}

\renewcommand{\thefootnote}{\arabic{footnote}}

\maketitle \centerline{December 18, 2009. Revised March 16, 2011} 

\begin{abstract}
We present in this paper two different classes of general $K$-splitting algorithms for
solving finite-dimensional convex optimization problems. Under the assumption that the function being minimized has a Lipschitz continuous gradient, we prove that the number of
iterations needed by the first class of algorithms to obtain an
$\epsilon$-optimal solution is $O(1/\epsilon)$. The algorithms in
the second class are accelerated versions of those in the first
class, where the complexity result is improved to
$O(1/\sqrt{\epsilon})$ while the computational effort required at
each iteration is almost unchanged. To the best
of our knowledge, the complexity results presented in this paper are
the first ones of this type that have been given for splitting and alternating
direction type methods. Moreover, all algorithms proposed in this
paper are parallelizable, which makes them particularly attractive
for solving certain large-scale problems.
\end{abstract}



\section{Introduction}\label{sec:intro} Many convex optimization problems that arise in practice take the form of a sum of convex functions.
Often one function is an energy that one wants to minimize and the
other functions are regularization terms to make the solution have certain properties. For example, Tikhonov regularization
\cite{Tikhonov-Arsenin-book-77} is usually applied to
ill-conditioned inverse problems to make them well-posed,
compressed sensing \cite{Candes-Romberg-Tao-2006,Donoho-06} uses
$\ell_1$ regularization to obtain sparse solutions, and problems arising from medical imaging adopt both $\ell_1$ and
total variation (TV) as regularization terms \cite{Ma-Yin-Zhang-Chakraborty-08}. In this paper,
we propose and analyze splitting/alternating direction
algorithms for solving the following convex optimization problem:
\bea\label{prob:convex-K-sum-original}\min_{x\in\br^n}\quad
F(x)\equiv\sum_{i=1}^K f_i(x),\eea where
$f_i:\br^n\rightarrow\br,i=1,\ldots,K,$ are convex functions. When the
functions $f_i$'s are well-structured, a well established way to
solve problem \eqref{prob:convex-K-sum-original} is to split the
variable $x$ into $K$ variables by introducing $K-1$ new variables
and then apply an augmented Lagrangian method to solve the resulting
problem. Decomposition of the augmented Lagrangian
function can then be accomplished by applying an alternating direction method (ADM) to minimize it.

Problem \eqref{prob:convex-K-sum-original} is closely related to the following inclusion problem:
\bea\label{prob:K-operator-zeros}0\in T_1(x)+\cdots+T_K(x),\eea where $T_1,\ldots,T_K$ are set-valued maximal monotone operators.
The goal of problem \eqref{prob:K-operator-zeros} is to find a zero of the sum of $K$ maximal monotone operators. Note that the optimality
conditions for \eqref{prob:convex-K-sum-original} are $$0\in \sum_{i=1}^K\partial f_i(x);$$ hence, these conditions can be satisfied by solving a problem
of the form \eqref{prob:K-operator-zeros}.

In the extensive literature on
splitting and ADM algorithms, the case $K=2$ predominates. The algorithms for solving \eqref{prob:K-operator-zeros} when $K=2$ are usually based on operator
splitting techniques. Important operator splitting algorithms include the Douglas-Rachford
\cite{Douglas-Rachford-56,Eckstein-Bertsekas-1992,Combettes-Pesquet-DR-2007}, Peaceman-Rachford \cite{Peaceman-Rachford-55},
double-backward \cite{Combettes-2004} and forward-backward class \cite{Gabay-83,Tseng-00} of algorithms. Alternating direction methods (ADM) within an
augmented Lagrangian framework for solving \eqref{prob:convex-K-sum-original} are optimization analogs/variants of the Douglas-Rachford and Peaceman-Rachford
splitting methods. These algorithms have been studied extensively for the case of $K=2$, and were first proposed in the 1970s for solving optimization problems
arising from numerical PDE problems \cite{Gabay-Mercier-1976,Glowinski-LeTallec-89}. We refer to \cite{Goldfarb-Ma-Scheinberg-2010} and the references therein for
more information on splitting and ADM algorithms for the case of $K=2$.

Although there is an extensive literature on operator splitting
methods, very few convergence results have been published on methods
for finding a zero of a sum of more than two maximal monotone
operators. The principal exceptions, are the Jacobi-like method of
Spingarn \cite{Spingarn-1983} and more recently, the general
projective splitting methods of Eckstein and Svaiter
\cite{Eckstein-Svaiter-2009}. The algorithm addressed in \cite{Spingarn-1983} first reduces problem \eqref{prob:K-operator-zeros} to the sum of two
maximal monotone operators by defining new subspaces and operators, and then applies a Douglas-Rachford splitting algorithm to solve the new problem.
The projective splitting methods in \cite{Eckstein-Svaiter-2009} do not reduce problem \eqref{prob:K-operator-zeros} to the case $K=2$.
Instead, by using the concept of an extended solution set, it is shown in \cite{Eckstein-Svaiter-2009} that solving \eqref{prob:K-operator-zeros} is equivalent
to finding a point in the extended solution set, and a separator-projection algorithm is given to do this.

Global convergence results for variable splitting ADMs and
operator splitting algorithms for the case of $K=2$ have been proved under
various assumptions. However, except for the fairly recently proposed gradient methods in \cite{Nesterov-07} and related iterative shrinkage/thresholding algorithms in \cite{Beck-Teboulle-2009} and the
alternating linearization methods in \cite{Goldfarb-Ma-Scheinberg-2010},
complexity bounds for these methods had not been established.
These complexity results are extensions of the seminal results of Nesterov \cite{Nesterov-1983,NesterovConvexBook2004},
who first showed that certain first-order methods that he proposed could obtain an $\epsilon$-optimal solution of a smooth convex programming problem
in $O(1/\sqrt{\epsilon})$ iterations. Moreover, he showed that his methods were optimal in the sense that this iteration complexity was the best that could be  obtained using only first-order information.
Nesterov's optimal gradient methods are
accelerated gradient methods that use a combination of previous points to compute the new point at each iteration. By combining
these methods with smoothing techniques, optimal complexity results
were obtained for solving nonsmooth problems in
\cite{Nesterov-2005,Tseng-2008}.

In this paper, we propose two classes of multiple variable-splitting algorithms based on alternating direction and alternating linearization techniques
that can solve problem \eqref{prob:convex-K-sum-original} for
general $K(K\geq 2)$ and we present complexity results for them. (Note that the complexity results in \cite{Nesterov-07,Beck-Teboulle-2009,Goldfarb-Ma-Scheinberg-2010} are only for problem \eqref{prob:convex-K-sum-original} when $K=2$). The algorithms in the
first class can be viewed as
alternating linearization methods in the sense that at each
iteration these algorithms perform $K$ minimizations of an approximation to the
original objective function $F$ by keeping one of the functions
$f_i(x)$ unchanged and linearizing the other $K-1$ functions. An
alternating linearization method for minimizing the sum of two
convex functions was studied by Kiwiel \etal
\cite{Kiwiel-Rosa-1999}. However, our algorithms differ greatly
from the one in \cite{Kiwiel-Rosa-1999} in the way that the proximal terms are chosen. Moreover, our
algorithms are more general as they can solve general problems with
$K(K\geq 2)$ functions. Furthermore, we prove that the iteration
complexity of this class of splitting algorithms is
$O(1/\epsilon)$ for an $\epsilon$-optimal solution. To the best of
our knowledge, this is the first complexity result of this type for
splitting/alternating direction type algorithms.
The algorithms in our second class are accelerated versions of the
algorithms in our first class and have $O(1/\sqrt{\epsilon})$ iteration
complexities. This class of splitting algorithms is also new as are
the complexity results.

Our new algorithms have, in addition, several practical advantages.
First, they are all parallelizable. Thus, although at each iteration
we solve $K$ subproblems, the CPU time required should be
approximately equal to the time required to solve the most difficult
of the subproblems if we have $K$ processors that can work in
parallel. Second, since every function $f_i$ is minimized once at each
iteration, it is likely that our algorithms will need fewer iterations
to converge than operator splitting algorithms such as FPC
\cite{Hale-Yin-Zhang-SIAM-2008,Ma-Goldfarb-Chen-2008},TVCMRI
\cite{Ma-Yin-Zhang-Chakraborty-08}, ISTA and FISTA
\cite{Beck-Teboulle-2009}. The numerical results in
\cite{Afonso-BD-Figueiredo-2009} for the case of $K=2$ support this conclusion.

The rest of this paper is organized as follows. In Section \ref{sec:MSA} we
propose a class of splitting algorithms based on alternating direction and alternating linearization methods for solving
\eqref{prob:convex-K-sum-original} and prove that they require $O(1/\epsilon)$ iterations to obtain an $\epsilon$-optimal solution. In
Section \ref{sec:FaMSA} we propose accelerated splitting algorithms for
solving \eqref{prob:convex-K-sum-original} and prove they have $O(1/\sqrt{\epsilon})$ complexities. We discuss how to
apply our algorithms for solving nonsmooth problems by using smoothing
techniques in Section \ref{sec:nonsmooth}. Numerical results are presented in Section \ref{sec:numerical}. Finally, we summarize our results in Section \ref{sec:discussion}.

\section{A class of multiple splitting algorithms}\label{sec:MSA} By introducing new variables, i.e., splitting variable $x$ into $K$ different variables,
problem \eqref{prob:convex-K-sum-original} can be rewritten as:
\bea\label{prob:convex-K-sum}\ba{ll}\min & \displaystyle\sum_{i=1}^K f_i(x^i) \\
s.t. & x^i=x^{i+1}, i=1,\ldots,K-1. \ea\eea In Sections 2 and 3, we
focus on splitting and ADM algorithms for solving
\eqref{prob:convex-K-sum} and their complexity results.

We make the following assumptions throughout Sections 2 and 3.
\begin{assumption}\label{assumption-section-2}
\begin{itemize}
\item $f_i(\cdot):\br^n\rightarrow\br,i=1,\ldots,K$ is a smooth convex function of the type $C^{1,1}$, i.e. continuously differentiable with
Lipschitz continuous gradient:
\begin{align*}\|\nabla f_i(x)-\nabla
f_i(y)\|\leq L(f_i)\|x-y\|, \forall x,y\in\br^n,\end{align*} where
$L(f_i)$ is the Lipschitz constant.
\item Problem \eqref{prob:convex-K-sum-original} is solvable, i.e., $X_*:=\arg\min F\neq\emptyset.$
\end{itemize}
\end{assumption}

We define the term $\epsilon$-optimal as follows.
\begin{definition}
Suppose $x^*$ is an optimal solution to the following problem
\bea\label{prob:epsilon-optimal}\min\{f(x):x\in \mathcal{C}\}.\eea $x\in \mathcal{C}$ is
called an $\epsilon$-optimal solution to
\eqref{prob:epsilon-optimal} if $f(x)-f(x^*)\leq \epsilon$ holds.
\end{definition}

The following notation is adopted throughout Sections 2 and 3.
\begin{definition}\label{def:Multiple-split}
We define $\tilde{f}_i(u,v)$
as the linear approximation to $f_i(u)$ at a point $v$ plus a
proximal term:
\begin{align*}\tilde{f}_i(u,v):=f_i(v)+\langle\nabla f_i(v),u-v\rangle+\oneotwomu\|u-v\|^2,\end{align*} where $\mu$ is a penalty parameter.
We use $Q_i(v^1,\ldots,v^{i-1},u,v^{i+1},\ldots,v^K)$ to denote the
following approximation to the function $F(u)$:
\begin{align*}Q_i(v^1,\ldots,v^{i-1},u,v^{i+1},\ldots,v^K):=f_i(u)+\sum_{j=1,j\neq i}^K \tilde{f}_j(u,v^j),
\end{align*} i.e., $Q_i$ is an approximation to the function $F$, where
the $i$-th function $f_i$ is unchanged but the other functions are
approximated by a linear term plus a proximal term. We use
$p_i(v^1,\ldots,v^{i-1},v^{i+1},\ldots,v^K)$ to denote the minimizer
of $Q_i(v^1,\ldots,v^{i-1},u,v^{i+1},\ldots,v^K)$ with respect to
$u$, i.e.,
\begin{align}\label{def:Multiple-split-def-pi}p_i(v^1,\ldots,v^{i-1},v^{i+1},\ldots,v^K):=\arg\min_u Q_i(v^1,\ldots,v^{i-1},u,v^{i+1},\ldots,v^K). \end{align}
\end{definition}

With the above notation, we have the following lemma which follows
from a fundamental property of a smooth function in the class
$C^{1,1}$; see e.g., \cite{Bertsekas-book-99}.
\begin{lemma}\label{lemma:Lipschtz-property}
For $\tilde{f}_i$ defined as in Definition \ref{def:Multiple-split} and $\mu\leq 1/\max_{1\leq i\leq K}{L(f_i)}$,
we have for $i=1,\ldots,K$,
\begin{align*} f_i(x)\leq f_i(y)+\langle\nabla f_i(y),x-y\rangle+\frac{L(f_i)}{2}\|x-y\|^2 \leq \tilde{f}_i(x,y), \forall x,y\in\br^n.\end{align*}
\end{lemma}

The following key lemma is crucial
for the proofs of our complexity results. Our proofs of this lemma
and most of the results that follow in this and the
remaining sections of the paper closely follow proofs given in
\cite{Beck-Teboulle-2009} for related lemmas and theorems.
\begin{lemma}\label{lemma:basic-multiple-split}
For any $i=1,\ldots,K$,
$u,v^1,\ldots,v^{i-1},v^{i+1},\ldots,v^K\in\br^n$ and $\mu\leq 1/\max_{1\leq i\leq K}{L(f_i)}$, we have,
\begin{align}\label{lemma:basic-multiple-split-inequa}2\mu(F(u)-F(p))\geq \sum_{j=1,j\neq i}^K \left(\|p-u\|^2-\|v^j-u\|^2\right), \end{align}
where $p:=p_i(v^1,\ldots,v^{i-1},v^{i+1},\ldots,v^K)$.
\end{lemma}
\begin{proof}
From Lemma \ref{lemma:Lipschtz-property} we know that $F(p)\leq
Q_i(v^1,\ldots,v^{i-1},p,v^{i+1},\ldots,v^K)$ holds for all $i$ and
$v^1,\ldots,v^{i-1},v^{i+1},\ldots,v^K\in\br^n$. Thus, for any
$u\in\br^n$ we have,
\begin{align}\label{lemma:basic-multiple-split-proof-1}F(u)-F(p) & \geq F(u)-Q_i(v^1,\ldots,v^{i-1},p,v^{i+1},\ldots,v^K) \\
                       \nonumber & = f_i(u)-f_i(p) + \sum_{j=1,j\neq i}^K \left(f_j(u)-f_j(v^j)-\left\langle\nabla
                        f_j(v^j),p-v^j\right\rangle+\oneotwomu\|p-v^j\|^2\right)\\
                       \nonumber & \geq \left\langle\nabla f_i(p),u-p\right\rangle + \sum_{j=1,j\neq
                        i}^K \left(\left\langle\nabla f_j(v^j),u-v^j\right\rangle-\left\langle\nabla
                        f_j(v^j),p-v^j\right\rangle+\oneotwomu\|p-v^j\|^2\right)\\
                       \nonumber & = \left\langle\nabla f_i(p),u-p\right\rangle + \sum_{j=1,j\neq
                        i}^K \left(\left\langle\nabla
                        f_j(v^j),u-p\right\rangle+\oneotwomu\|p-v^j\|^2\right)\\
                      \nonumber  & = \sum_{j=1,j\neq i}^K
                        \left(\left\langle-\oneomu(p-v^j),u-p\right\rangle-\oneotwomu\|p-v^j\|^2\right),
\end{align}
where the second inequality is due to the convexity of the functions
$f_j,j=1,\ldots,K$ and the last equality is from the first-order
optimality conditions for problem \eqref{def:Multiple-split-def-pi},
i.e., \bea\label{def:Multiple-split-def-pi-optcond}\nabla f_i(p)+\sum_{j=1,j\neq i}^K\left(\nabla
f_j(v^j)+\oneomu(p-v^j)\right)=0.\eea Then using the identity
\bea\label{Pythagoras-identity}\|a-c\|^2-\|b-c\|^2=\|a-b\|^2+2\langle a-b,b-c\rangle,\eea
we get the following
inequality:
\begin{align*}2\mu(F(u)-F(p)) & \geq \sum_{j=1,j\neq i}^K
\left(\left\langle-2(p-v^j),u-p\right\rangle-\|p-v^j\|^2\right)\\
& = \sum_{j=1,j\neq i}^K \left(\|p-u\|^2-\|v^j-u\|^2\right).
\end{align*}
\end{proof}

Our multiple splitting algorithms (MSA) for
solving \eqref{prob:convex-K-sum} are outlined in
Algorithm \ref{alg:multiple-splitting}, where $D^{(k)}\in\br^{K\times
K}$ is a doubly stochastic matrix, i.e.,
$$D^{(k)}_{ij}\geq 0, \sum_{j=1}^K D^{(k)}_{ij}=1,
\sum_{i=1}^K D^{(k)}_{ij}=1,\forall i,j=1,\ldots,K.$$
\begin{algorithm2e}\caption{A Class of Multiple Splitting Algorithms (MSA)}
    \label{alg:multiple-splitting}
\linesnumberedhidden \dontprintsemicolon Set
$x_{0}=x_{(0)}^1=\ldots=x_{(0)}^K=w_{(0)}^1=\ldots=w_{(0)}^K$ and $\mu\leq 1/\max_{1\leq i\leq K}{L(f_i)}$.\;
\For{$k=0,1,\cdots$}{
\begin{itemize}
\item for each $i=1,\ldots,K$, compute
$x_{(k+1)}^i := p_i(w_{(k)}^i,\ldots,w_{(k)}^{i}),$\; 
\item compute
\begin{align*}
\begin{pmatrix}w_{(k+1)}^1, \ldots,
w_{(k+1)}^K\end{pmatrix} & := \begin{pmatrix} x_{(k+1)}^1,
\ldots,x_{(k+1)}^K \end{pmatrix}D^{(k+1)}.
\end{align*}
\end{itemize}}
\end{algorithm2e}
One natural
choice of $D^{(k)}$ is to take all of its components equal to $1/K$.
In this case, all $w_{(k)}^i,i=1,\ldots,K$ are equal to
$\sum_{i=1}^K x_{(k)}^i/K$, i.e., the average of the current $K$
iterates.

At iteration $k$, Algorithm \ref{alg:multiple-splitting} computes
$K$ points $x_{(k)}^i,i=1,\ldots,K$ by solving $K$ subproblems. For
many problems in practice, these $K$ subproblems are expected to be
very easy to solve. Another advantage of the algorithm is that it is
parallelizable since given $w_{(k)}^i,i=1,\ldots,K$, the $K$
subproblems in Algorithm \ref{alg:multiple-splitting} can be solved
simultaneously. Algorithm \eqref{alg:multiple-splitting} can be
viewed as an alternating linearization method since at each
iteration, $K$ subproblems are solved, and each subproblem
corresponds to minimizing a function involving linear approximations
to some of the functions. Although Algorithm
\ref{alg:multiple-splitting} assumes the Lipschitz constants are
known, and hence that $\mu$ is known, this assumption can be relaxed by using
the backtracking technique in \cite{Beck-Teboulle-2009} to estimate
$\mu$ at each iteration.

We prove in the following that the number of
iterations needed by Algorithm \ref{alg:multiple-splitting} to obtain an $\epsilon$-optimal solution is
$O(1/\epsilon)$.

\begin{theorem}\label{the:multiple-splitting}
Suppose $x^*$ is an optimal solution to problem
\eqref{prob:convex-K-sum}. For any choice of $\mu\leq 1/\max_{1\leq i\leq K}{L(f_i)}$, the sequence
$\{x_{(k)}^i,w_{(k)}^i\}_{i=1}^K$ generated by Algorithm
\ref{alg:multiple-splitting} satisfies:
\bea\label{the:multiple-splitting-conclusion}\min_{i=1,\ldots,K}F(x_{(k)}^i)-F(x^*)\leq
\frac{(K-1)\|x_{0}-x^*\|^2}{2\mu k}.\eea Thus, the sequence $\{\min_{i=1,\ldots,K}F(x_{(k)}^i)\}$ produced by
Algorithm \ref{alg:multiple-splitting} converges to $F(x^*)$. Moreover, if $\mu\geq\beta/\max_i\{L(f_i)\}$ where $0<\beta\leq 1$,
the number of iterations needed to obtain an $\epsilon$-optimal solution is at most $\lceil C/\epsilon \rceil$,
where $C=\frac{(K-1)\max_i\{L(f_i)\}\|x_0-x^*\|^2}{2\beta}$.
\end{theorem}

\begin{proof}
In \eqref{lemma:basic-multiple-split-inequa}, by letting
$u=x^*,v^j=w_{(n)}^i,j=1,\ldots,K,j\neq i$, we have $p=x_{n+1}^i$
and
\begin{align}\label{proof:multiple-splitting-ineq-1}2\mu(F(x^*)-F(x_{(n+1)}^i))&\geq
\sum_{j=1,j\neq
i}^K\left(\|x_{(n+1)}^i-x^*\|^2-\|w_{(n)}^i-x^*\|^2\right)\\\nonumber
& =
(K-1)\left(\|x_{(n+1)}^i-x^*\|^2-\|w_{(n)}^i-x^*\|^2\right).\end{align}
Using the definition of $w_{(n)}^i$ in Algorithm
\ref{alg:multiple-splitting}, we have
\begin{align}\label{proof:multiple-splitting-ineq-2}\sum_{i=1}^K\|w_{(n+1)}^i-x^*\|^2 & = \sum_{i=1}^K\left\|\sum_{j=1}^K D_{ji}^{(n+1)}x_{(n+1)}^j-x^*\right\|^2
           \\\nonumber  & = \sum_{i=1}^K\left\|\sum_{j=1}^K D_{ji}^{(n+1)}(x_{(n+1)}^j-x^*)\right\|^2
            \\\nonumber & \leq \sum_{i=1}^K\sum_{j=1}^K
            D_{ji}^{(n+1)}\|x_{(n+1)}^j-x^*\|^2 \\\nonumber
            & = \sum_{j=1}^K\|x_{(n+1)}^j-x^*\|^2,
           \end{align}
where the second and the last equalities are due to the fact that
$D^{(n+1)}$ is a doubly stochastic matrix and the inequality is due
to the convexity of the function $\|\cdot\|^2.$

Thus by summing \eqref{proof:multiple-splitting-ineq-1} over
$i=1,\ldots,K$ we obtain
\begin{align}\label{proof:multiple-splitting-ineq-3}2\mu\left(KF(x^*)-\sum_{i=1}^KF(x_{(n+1)}^i)\right)&
\geq(K-1)\left(\sum_{i=1}^K\|x_{(n+1)}^i-x^*\|^2-\sum_{i=1}^K\|w_{(n)}^i-x^*\|^2\right)
\\\nonumber & \geq (K-1)\left(\sum_{i=1}^K\|w_{(n+1)}^i-x^*\|^2-\sum_{i=1}^K\|w_{(n)}^i-x^*\|^2\right),
\end{align}
where the last inequality is due to
\eqref{proof:multiple-splitting-ineq-2}.

Summing \eqref{proof:multiple-splitting-ineq-3} over
$n=0,1,\ldots,k-1$, and using the fact that
$w_{(0)}^i=x_{0},i=1,\ldots,K$, yields
\begin{align}\label{proof:multiple-splitting-ineq-4}2\mu\left(kKF(x^*)-\sum_{n=0}^{k-1}\sum_{i=1}^KF(x_{(n+1)}^i)\right)
& \geq
(K-1)\sum_{i=1}^K\left(\|w_{(k)}^i-x^*\|^2-\|w_{(0)}^i-x^*\|^2\right)
\\ \nonumber & \geq -K(K-1)\|x_{0}-x^*\|^2.\end{align}

In \eqref{lemma:basic-multiple-split-inequa}, by letting
$u=v^j=w_{(n)}^i, j=1,\ldots,K,j\neq i$, we get $p=x_{(n+1)}^i$ and
\bea\label{proof:multiple-splitting-ineq-5}2\mu\left(F(w_{(n)}^i)-F(x_{(n+1)}^i)\right)\geq
\sum_{j=1,j\neq i}^K\|x_{(n+1)}^i-w_{(n)}^i\|^2\geq 0.\eea From the
way we compute $w_{(n)}^i, i=1,\ldots,K,$ and the facts that $F$ is
convex and $D^{(n)}$ is a doubly stochastic matrix, we get
\begin{align}\label{proof:multiple-splitting-ineq-6}\sum_{i=1}^KF(w_{(n)}^i) & = \sum_{i=1}^KF\left(\sum_{j=1}^KD_{ji}^{(n)}x_{(n)}^j\right)\\
             \nonumber & \leq \sum_{i=1}^K\sum_{j=1}^KD_{ji}^{(n)}F(x_{(n)}^j) \\
             \nonumber & = \sum_{j=1}^KF(x_{(n)}^j).\end{align}

Now summing \eqref{proof:multiple-splitting-ineq-5} over
$i=1,\ldots,K$ and using \eqref{proof:multiple-splitting-ineq-6}, we
obtain
\begin{align}\label{proof:multiple-splitting-ineq-7}2\mu\left(\sum_{i=1}^KF(x_{(n)}^i)-\sum_{i=1}^KF(x_{(n+1)}^i)\right)\geq 0.\end{align}
This shows that the sums $\sum_{i=1}^KF(x_{(n)}^i)$ are non-increasing as $n$ increases. Hence,
\begin{align}\label{proof:multiple-splitting-ineq-8} \sum_{n=0}^{k-1}\sum_{i=1}^KF(x_{(n+1)}^i)\geq k\sum_{i=1}^KF(x_{(k)}^i).\end{align}

Finally, combining \eqref{proof:multiple-splitting-ineq-4} and
\eqref{proof:multiple-splitting-ineq-8} yields
\begin{align*}2\mu\left(kKF(x^*)-k\sum_{i=1}^KF(x_{(k)}^i)\right)\geq-K(K-1)\|x_{0}-x^*\|^2.\end{align*}
Hence,
\begin{align*}\min_{1,\ldots,K}F(x_{(k)}^i)-F(x^*)\leq\frac{1}{K}\sum_{i=1}^KF(x_{(k)}^i)-F(x^*)\leq \frac{(K-1)\|x_{0}-x^*\|^2}{2\mu k}.\end{align*}

It then follows that $\min_{i=1,\ldots,K} F(x_{(k)}^{i})-F(x^*)\leq
\epsilon,$ if $k\geq \lceil C/\epsilon\rceil$, where
$C=\frac{(K-1)\max_i\{L(f_i)\}\|x_0-x^*\|^2}{2\beta}$, and hence that for any $k\geq \lceil C/\epsilon\rceil$, $x_{(k)}:=\arg\min\{x_{(k)}^i | F(x_{(k)}^i),
i=1,\ldots,K\}$ is an $\epsilon$-optimal solution.
\end{proof}

\begin{remark}
If in the original problem \eqref{prob:convex-K-sum-original}, $x$ is subject to a convex constraint $x\in\mathcal{C}$, where $\mathcal{C}$ is a convex set,
we can impose this constraint in every subproblem in MSA and obtain the same complexity result. The only changes in the proof are in
Lemma \ref{lemma:basic-multiple-split}. If there is a constraint $x\in\mathcal{C}$, then \eqref{lemma:basic-multiple-split-inequa} and
\eqref{lemma:basic-multiple-split-proof-1} hold for any $u\in\mathcal{C}$ and the last equality in \eqref{lemma:basic-multiple-split-proof-1}
becomes a ``$\geq$'' inequality due to the fact that the optimality conditions \eqref{def:Multiple-split-def-pi-optcond} become
\beaa \left\langle \nabla f_i(p)+\sum_{j=1,j\neq i}^K\left(\nabla
f_j(v^j)+\oneomu(p-v^j)\right), u-p\right\rangle \geq  0, \forall u\in\mathcal{C}.\eeaa Unfortunately, this extension is not very practical,
since for it to be useful, adding the constraint in every subproblem would most likely make most of these subproblems difficult to solve.
\end{remark}

\section{A class of fast multiple splitting algorithms}\label{sec:FaMSA} In this section, we give a class of fast multiple splitting
algorithms (FaMSA) for solving problem \eqref{prob:convex-K-sum} that
require at most $O(1/\sqrt{\epsilon})$ iterations to obtain an
$\epsilon$-optimal solution while requiring a computational effort
at each iteration that is roughly the same as Algorithm
\ref{alg:multiple-splitting}. Our
fast multiple splitting algorithms are outlined in Algorithm
\ref{alg:multiple-splitting-faster}, where $D^{(k)}\in\br^{K\times K}$ is a doubly stochastic matrix.
\begin{algorithm2e}\caption{A Class of Fast Multiple Splitting Algorithms (FaMSA)}
    \label{alg:multiple-splitting-faster}
\linesnumberedhidden \dontprintsemicolon Set
$x_{0}=x_{(0)}^i=\hat{w}_{(0)}^i=w_{(1)}^i, i=1,\ldots,K, t_1=1,$ and choose $\mu\leq 1/\max_{1\leq i\leq K}{L(f_i)}$ \;
\For{$k=1,2,\cdots$}{
\begin{itemize}
\item for $i=1,\ldots,K$, compute
$x_{(k)}^i=p_i(w_{(k)}^i,\ldots,w_{(k)}^i)$
\item compute $\begin{pmatrix}\hat{w}_{(k)}^1,\ldots,\hat{w}_{(k)}^K\end{pmatrix}:=\begin{pmatrix}x_{(k)}^1,\ldots,x_{(k)}^K\end{pmatrix}D^{(k)}$
\item compute $t_{k+1}=(1+\sqrt{1+4t_k^2})/2$
\item for $i=1,\ldots,K$, compute $w_{(k+1)}^i:=\hat{w}_{(k)}^i+\frac{1}{t_{k+1}}\left(t_k(x_{(k)}^i-\hat{w}_{(k-1)}^i)-(\hat{w}_{(k)}^i-\hat{w}_{(k-1)}^i)\right).$
\end{itemize}
}
\end{algorithm2e}

To establish the $O(1/\sqrt{\epsilon})$ iteration complexity of FaMSA, we need the
following lemma.

\begin{lemma}\label{lemma:multiple-splitting-faster}
Suppose $x^*$ is an optimal solution to problem
\eqref{prob:convex-K-sum}. For any choice of $\mu\leq\max_{1\leq i\leq K} L(f_i)$, the sequence $\{x_{(k)}^i,w_{(k)}^i,
\hat{w}_{(k)}^i\}_{i=1}^K$ generated by Algorithm
\ref{alg:multiple-splitting-faster} satisfies:
\bea\label{lemma-multiple-splitting-faster-inequa}2\mu(t_k^2v_k-t_{k+1}^2v_{k+1})\geq
(K-1)\sum_{i=1}^K \left(\|u_{k+1}^i\|^2-\|u_k^i\|^2\right),\eea
where $v_k:=\sum_{i=1}^K F(x_{(k)}^i)-KF(x^*)$ and
$u_k^i:=t_kx_{(k)}^i-(t_k-1)\hat{w}_{(k-1)}^i-x^*, i=1,\ldots,K.$
\end{lemma}

\begin{proof}
In \eqref{lemma:basic-multiple-split-inequa}, by letting
$u=\hat{w}_{(k)}^i,v^j=w_{(k+1)}^i,j=1,\ldots,K,j\neq i$, we get
$p=x_{(k+1)}^i$ and
\begin{align}\label{proof:multiple-splitting-faster-lemma-1}2\mu\left(F(\hat{w}_{(k)}^i)-F(x_{(k+1)}^i)\right)&\geq
\sum_{j=1,j\neq
i}^K\left(\|x_{(k+1)}^i-\hat{w}_{(k)}^i\|^2-\|w_{(k+1)}^i-\hat{w}_{(k)}^i\|^2\right)\\\nonumber
& =
(K-1)\left(\|x_{(k+1)}^i-\hat{w}_{(k)}^i\|^2-\|w_{(k+1)}^i-\hat{w}_{(k)}^i\|^2\right).\end{align}

Summing \eqref{proof:multiple-splitting-faster-lemma-1} over
$i=1,\ldots,K$, and using the facts that $F$ is convex and $D^{(k)}$ is a doubly stochastic matrix, we get
\begin{align*}2\mu\left(\sum_{i=1}^KF(x_{(k)}^i)-\sum_{i=1}^KF(x_{(k+1)}^i)\right)&\geq 2\mu\left(\sum_{i=1}^KF(\hat{w}_{(k)}^i)-\sum_{i=1}^KF(x_{(k+1)}^i)\right)
\\\nonumber & \geq (K-1)\sum_{i=1}^K\left(\|x_{(k+1)}^i-\hat{w}_{(k)}^i\|^2-\|w_{(k+1)}^i-\hat{w}_{(k)}^i\|^2\right), \end{align*}
i.e.,
\begin{align}\label{proof:multiple-splitting-faster-lemma-2} 2\mu(v_k-v_{k+1}) \geq (K-1)\sum_{i=1}^K
\left(\|x_{(k+1)}^i-\hat{w}_{(k)}^i\|^2-\|w_{(k+1)}^i-\hat{w}_{(k)}^i\|^2\right).\end{align}

In \eqref{lemma:basic-multiple-split-inequa}, by letting
$u=x^*,v^j=w_{(k+1)}^i$, we get $p=x_{(k+1)}^i$ and
\begin{align}\label{proof:multiple-splitting-faster-lemma-3}2\mu\left(F(x^*)-F(x_{(k+1)}^i)\right)&\geq \sum_{j=1,j\neq i}^K
\left(\|x_{(k+1)}^i-x^*\|^2-\|w_{(k+1)}^i-x^*\|^2\right)
\\\nonumber & = (K-1)\left(\|x_{(k+1)}^i-x^*\|^2-\|w_{(k+1)}^i-x^*\|^2\right).\end{align}
Summing \eqref{proof:multiple-splitting-faster-lemma-3} over
$i=1,\ldots,K$ we obtain
\begin{align}\label{proof:multiple-splitting-faster-lemma-4}-2\mu
v_{k+1}\geq
(K-1)\sum_{i=1}^K\left(\|x_{(k+1)}^i-x^*\|^2-\|w_{(k+1)}^i-x^*\|^2\right).\end{align}

Now multiplying \eqref{proof:multiple-splitting-faster-lemma-2} by
$t_k^2$ and
\eqref{proof:multiple-splitting-faster-lemma-4} by $t_{k+1}$,
adding the resulting two inequalities, using the relation
$t_k^2=t_{k+1}(t_{k+1}-1)$, and the identity \eqref{Pythagoras-identity}, we get
\begin{align}\label{proof:multiple-splitting-faster-lemma-5}
 & 2\mu(t_k^2v_k-t_{k+1}^2v_{k+1}) \\ \nonumber \geq &
 (K-1)\sum_{i=1}^K
 t_{k+1}(t_{k+1}-1)\left(\|x_{(k+1)}^i-\hat{w}_{(k)}^i\|^2-\|w_{(k+1)}^i-\hat{w}_{(k)}^i\|^2\right)
 \\ \nonumber & + (K-1)\sum_{i=1}^K
 t_{k+1}\left(\|x_{(k+1)}^i-x^*\|^2-\|w_{(k+1)}^i-x^*\|^2\right)
 \\\nonumber = & (K-1)\sum_{i=1}^K t_{k+1}(t_{k+1}-1)\left(\|x_{(k+1)}^i-w_{(k+1)}^i\|^2+2\left
 \langle x_{(k+1)}^i-w_{(k+1)}^i,w_{(k+1)}^i-\hat{w}_{(k)}^i\right\rangle\right)
 \\\nonumber & +(K-1)\sum_{i=1}^K t_{k+1}\left(\|x_{(k+1)}^i-w_{(k+1)}^i\|^2+2\left\langle x_{(k+1)}^i-w_{(k+1)}^i,w_{(k+1)}^i-x^*\right\rangle \right)
 \\\nonumber
 = & (K-1)\sum_{i=1}^K \left(t_{k+1}^2\|x_{(k+1)}^i-w_{(k+1)}^i\|^2+2t_{k+1}\left
 \langle x_{(k+1)}^i-w_{(k+1)}^i,t_{k+1}w_{(k+1)}^i-(t_{k+1}-1)\hat{w}_{(k)}^i-x^* \right\rangle\right)
 \\ \nonumber = & (K-1)\sum_{i=1}^K \left(\|t_{k+1}x_{(k+1)}^i-(t_{k+1}-1)\hat{w}_{(k)}^i-x^*\|^2-\|t_{k+1}w_{(k+1)}^i-(t_{k+1}-1)\hat{w}_{(k)}^i-x^*\|^2\right).
\end{align}
From the way we compute $w_{(k+1)}^i$ in Algorithm \ref{alg:multiple-splitting-faster}, i.e.,
$$w_{(k+1)}^i:=\hat{w}_{(k)}^i+\frac{1}{t_{k+1}}\left(t_k(x_{(k)}^i-\hat{w}_{(k-1)}^i)-(\hat{w}_{(k)}^i-\hat{w}_{(k-1)}^i)\right),$$ it follows that
$$t_{k+1}w_{(k+1)}^i-(t_{k+1}-1)\hat{w}_{(k)}^i-x^*=t_{k}x_{(k)}^i-(t_{k}-1)\hat{w}_{(k-1)}^i-x^*.$$
Thus, from \eqref{proof:multiple-splitting-faster-lemma-5} and the definition of $u_k^i$ it follows that
\begin{align*}& 2\mu(t_k^2v_k-t_{k+1}^2v_{k+1}) \\ \nonumber \geq & (K-1)\sum_{i=1}^K
\left(\|t_{k+1}x_{(k+1)}^i-(t_{k+1}-1)\hat{w}_{(k)}^i-x^*\|^2-\|t_{k}x_{(k)}^i-(t_{k}-1)\hat{w}_{(k-1)}^i-x^*\|^2\right)
\\\nonumber = & (K-1)\sum_{i=1}^K \left(\|u_{k+1}^i\|^2-\|u_k^i\|^2\right). \end{align*} This
completes the proof.
\end{proof}

Before proving our main complexity theorem to Algorithm \ref{alg:multiple-splitting-faster}, we note that the
sequence $\{t_k\}$ generated by Algorithm \ref{alg:multiple-splitting-faster} clearly satisfies $t_{k+1}\geq t_k+\half,$
and hence $t_k\geq (k+1)/2$ for all $k\geq 1$ since $t_1=1$.
\begin{theorem}\label{the:multiple-splitting-fast}
Suppose $x^*$ is an optimal solution to problem
\eqref{prob:convex-K-sum}. For any choice of $\mu\leq\max_{1\leq i\leq K}L(f_i)$, the sequence
$\{x_{(k)}^i,w_{(k)}^i,\hat{w}_{(k)}^i\}_{i=1}^K$ generated by Algorithm
\ref{alg:multiple-splitting-faster} satisfies:
\begin{align}\label{the:multiple-splitting-faster-conclusion}\min_{i=1,\ldots,K}F(x_{(k)}^i)-F(x^*)\leq
\frac{2(K-1)\|x_0-x^*\|^2}{\mu(k+1)^2}.\end{align} Thus, the sequence $\{\min_{i=1,\ldots,K}F(x_{(k)}^i)\}$
produced by Algorithm \ref{alg:multiple-splitting-faster} converges to $F(x^*)$. Moreover, if $\mu\geq\beta/\max_i\{L(f_i)\}$
where $0<\beta\leq 1$, the number of iterations needed to obtain an $\epsilon$-optimal solution is at most $\lfloor \sqrt{C/\epsilon} \rfloor$,
where $C=2(K-1)\max_i\{L(f_i)\}\|x_0-x^*\|^2/\beta$.
\end{theorem}
\begin{proof}
By rewriting \eqref{lemma-multiple-splitting-faster-inequa} as \beaa 2\mu t_{k+1}^2v_{k+1}+(K-1)
\sum_{i=1}^K\|u_{k+1}^i\|^2 \leq 2\mu t_k^2v_k+(K-1)\sum_{i=1}^K\|u_k^i\|^2, \eeaa we get
\begin{align*} 2\mu \left(\frac{k+1}{2}\right)^2 v_k & \leq 2\mu t_k^2 v_k+(K-1)\sum_{i=1}^K\|u_k^i\|^2 \\ &
\leq 2\mu t_1^2 v_1 + (K-1)\sum_{i=1}^K\|u_1^i\|^2 \\ & = 2\mu v_1+(K-1)\sum_{i=1}^K\|x_{(1)}^i-x^*\|^2 \\ &
\leq (K-1)\sum_{i=1}^K\|w_{(1)}^i-x^*\|^2 \\ & = K(K-1)\|x_0-x^*\|^2,\end{align*}
where the first inequality is due to $t_k\geq (k+1)/2$, the first equality is from the facts that $t_1=1$ and
$u_1^i=x_{(1)}^i-x^*$, the third inequality is from letting $k=0$ in \eqref{proof:multiple-splitting-faster-lemma-4}
and the last equality is due to $w_{(1)}^i=x_0,i=1,\ldots,K.$

Thus, from $v_k=\sum_{i=1}^K F(x_{(k)}^i)-K F(x^*)$ we get \beaa \sum_{i=1}^K F(x_{(k)}^i)-K F(x^*)\leq\frac{2K(K-1)\|x_0-x^*\|^2}{\mu(k+1)^2},\eeaa
which implies that \beaa \min_{i=1,\ldots,K}F(x_{(k)}^i)-F(x^*)\leq \frac{1}{K}\sum_{i=1}^K F(x_{(k)}^i)-F(x^*)\leq\frac{2(K-1)\|x_0-x^*\|^2}{\mu(k+1)^2}, \eeaa
i.e., \eqref{the:multiple-splitting-faster-conclusion} holds.

Moreover, it follows that if $C/(k+1)^2\leq \epsilon$, i.e., $k\geq
\lfloor\sqrt{C/\epsilon}\rfloor$, then $\min_{i=1,\ldots,K} F(x_{(k)}^{i})-F(x^*)\leq \epsilon,$ where
$C=2(K-1)\max_i\{L(f_i)\}\|x_0-x^*\|^2/\beta$. This implies that for any
$k\geq \lfloor\sqrt{C/\epsilon}\rfloor$, $x_{(k)}:=\arg\min\{x_{(k)}^i |
F(x_{(k)}^i), i=1,\ldots,K\}$ is an $\epsilon$-optimal solution.
\end{proof}

\begin{remark}
Although we have assumed that the Lipschitz constants $L(f_i)$ are known, and hence that $\mu$ is chosen in Algorithm \ref{alg:multiple-splitting-faster} to be smaller than $1/\max_{1\leq i\leq K}\{L(f_i)\}$,
this can be relaxed by using the backtracking technique in \cite{Beck-Teboulle-2009} that chooses a $\mu$ at each iteration that is smaller than the $\mu$ used at the previous iteration and for which $F(p)\leq Q_i(w_{(k)}^i,\ldots,w_{(k)}^i,p,w_{(k)}^i,\ldots,w_{(k)}^i)$ for all $i$.
\end{remark}

\subsection{A variant of the fast multiple splitting algorithm}
In this section, we present a variant of the fast multiple splitting algorithm (Algorithm \ref{alg:multiple-splitting-faster}) that is much more efficient and requires much less memory than Algorithm \ref{alg:multiple-splitting-faster} for problems in which $K$ is large. This variant uses $D^{(k)}:=1/Kee^\top$, where $e$ is the $n$-dimensional vector with all ones, and replaces $x_{(k)}^i$ in the last line of Algorithm \ref{alg:multiple-splitting-faster} by $\hat{w}_{(k)}^i$; i.e., in the last line of Algorithm \ref{alg:multiple-splitting-faster}, we compute $w_{(k+1)}^i$ for $i=1,\ldots,K$ by the formula:
\[w_{(k+1)}^i:=\hat{w}_{(k)}^i+\frac{t_k-1}{t_{k+1}}(\hat{w}_{(k)}^i-\hat{w}_{(k-1)}^i).\] It is easy to see that in this variant, the $\hat{w}_{(k)}^i,i=1,\ldots,K$ are all the same and the $w_{(k+1)}^i,i=1,\ldots,K$ are all the same. We call this variant FaMSA-s, where s refers to the fact that this variant computes a ``single'' vector $\hat{w}^k$ and a single vector $w_{(k+1)}$ at the $k$-th iteration. It is given below as Algorithm \ref{alg:multiple-splitting-faster-v}.
\begin{algorithm2e}\caption{A variant of FaMSA (FaMSA-s)}
    \label{alg:multiple-splitting-faster-v}
\linesnumberedhidden \dontprintsemicolon Set
$x_{0}=x_{(0)}^i=\hat{w}_{(0)}=w_{(1)}, i=1,\ldots,K, t_1=1,$ and choose $\mu\leq 1/\max_{1\leq i\leq K}{L(f_i)}$ \;
\For{$k=1,2,\cdots$}{
\begin{itemize}
\item for $i=1,\ldots,K$, compute
$x_{(k)}^i=p_i(w_{(k)},\ldots,w_{(k)})$
\item compute $\hat{w}_{(k)}:=\frac{1}{K}\sum_{i=1}^K x_{(k)}^i$
\item compute $t_{k+1}=(1+\sqrt{1+4t_k^2})/2$
\item for $i=1,\ldots,K$, compute $w_{(k+1)}:=\hat{w}_{(k)}+\frac{t_k-1}{t_{k+1}}(\hat{w}_{(k)}-\hat{w}_{(k-1)}).$
\end{itemize}
}
\end{algorithm2e}

It is easy to verify that the following analog of Lemma \ref{lemma:multiple-splitting-faster} applies to Algorithm FaMSA-s.

\begin{lemma}\label{lemma:multiple-splitting-faster-variant}
Suppose $x^*$ is an optimal solution to problem
\eqref{prob:convex-K-sum}. The sequence $\{w_{(k)},
\hat{w}_{(k)}\}$ generated by Algorithm
FaMSA-s satisfies:
\beaa 2\mu(t_k^2v_k-t_{k+1}^2v_{k+1})\geq
(K-1)\sum_{i=1}^K \left(\|u_{k+1}\|^2-\|u_k\|^2\right),\eeaa
where $v_k:= K(F(\hat{w}_{(k)})-F(x^*))$ and
$u_k:=t_k\hat{w}_{(k)}-(t_k-1)\hat{w}_{(k-1)}-x^*.$
\end{lemma}

\begin{proof}
The proof is very similar to the proof of Lemma \ref{lemma:multiple-splitting-faster}; hence, we leave it to the reader. The main difference is that instead of using the inequality $\sum_{i=1}^K F(\hat{w}_{(k)}^i)\leq \sum_{i=1}^K F(x_{(k)}^i)$ to replace the sum involving $\hat{w}_{(k)}^i$, we use the fact that $KF(\hat{w}_{k+1})\leq\sum_{i=1}^KF(x_{(k+1)}^i)$ to replace the sum involving $x_{(k+1)}^i$ in the proof.
\end{proof}

From Lemma \ref{lemma:multiple-splitting-faster-variant}, Theorem \ref{the:multiple-splitting-fast} with $\hat{w}_{(k)}^i$ and $w_{(k)}^i$, respectively, for $i=1,\ldots,K$ replaced by $\hat{w}_{(k)}$ and $w_{(k)}$ follows immediately for FaMSA-s.

\section{Multiple splitting algorithms for nonsmooth problems}\label{sec:nonsmooth}
Although for the above results we required all functions to be
in the class of $C^{1,1}$, our algorithms can still be applied to
solve nonsmooth problems by first smoothing all nonsmooth functions.
One of the most important smoothing techniques is the one proposed
by Nesterov \cite{Nesterov-2005}. We use the $\ell_1$-norm function as an
example to show how Nesterov's smoothing technique works. Note that
the $\ell_1$ function $f(x):=\|x\|_1$ can be rewritten as
$\|x\|_1=\max\{\langle x,u\rangle: u\in U\},$ where
$U:=\{u:\|u\|_\infty\leq 1\}.$ Since $U$ is a bounded convex set, we
can define a prox-function $d(u)$ for the set $U$, where $d(u)$ is continuous and strongly convex on $U$ with
convexity parameter $\sigma>0.$ For $U$ defined as above, a natural
choice for $d(u)$ is $d(u):=\half\|u\|_2^2$ and thus $\sigma=1.$ Hence, we have the following smooth approximation for $f(x)=\|x\|_1$:
\begin{align*}f_{\rho}(x):=\max\{\langle x,u\rangle-\rho d(u): u\in U\},\end{align*}
where $\rho$ is a positive smoothness parameter. It can be shown
that $f_{\rho}(x)$ is well defined and is in the class of
$C^{1,1}$ and its gradient is Lipschitz continuous with constant
$L_{\rho}=\frac{1}{\rho\sigma}$ (see Theorem 1 in
\cite{Nesterov-2005}). Also, it is easy to show that the following
relations hold for $f(x)$ and $f_{\rho}(x)$:
\begin{align*}f_{\rho}(x)
\leq f(x)\leq f_{\rho}(x)+\rho D,\end{align*} where $D:=\max
\limits_u\{d(u):u\in U\}. $ Therefore, to get an $\epsilon$-optimal
solution to a problem involving the $\ell_1$-norm function $f(x)$, we can
replace $f(x)$ with $f_{\frac{\epsilon}{2D}}(x)$ to get a smooth problem, and
then apply our splitting algorithms to
the new problem to get an $\frac{\epsilon}{2}$-optimal solution, which will be
$\epsilon$-optimal to the original nonsmooth problem. Since $L_{\frac{\epsilon}{2D}}$ is $O(1/\epsilon)$, our fast $O(1/\sqrt{\epsilon})$
algorithms require $O(1/\epsilon)$ iterations to compute an $\epsilon$-optimal solution.

For nonsmooth problems in imaging, data analysis, and machine learning, etc. with regularization terms that involve total variation and the
nuclear norm, we can use similar smoothing techniques to smooth these nonsmooth functions, and then apply our multiple splitting algorithms to solve them.

\section{Numerical experiments}\label{sec:numerical}
We present some preliminary numerical experiments in this section. Specifically, we apply our MSA and FaMSA algorithms to solve the Fermat-Weber
problem and a total variation and wavelet based image deblurring problem. All numerical experiments were run in MATLAB 7.3.0 on a Dell Precision 670 workstation with an Intel Xeon(TM) 3.4GHZ CPU and 6GB of RAM.

\subsection{The Fermat-Weber problem}
The Fermat-Weber (F-W) problem can be cast as:
\bea\label{prob:weber}\min F(x)\equiv\sum_{i=1}^K\|x-c^i\|,\eea where $c^i\in\br^n, i=1,\ldots,K$ are $K$ given points. Problem \eqref{prob:weber}
can be reformulated as a second-order cone programming (SOCP) problem and thus solved in polynomial time by an interior-point method. Since there are $K$ cones, the size of a standard form SOCP formulation for this problem is quite large for large $K$ and $n$.
Since $f_i(x)=\|x-c^i\|,i=1,\ldots,K$ are not smooth, to apply our MSA and FaMSA algorithms, we need to smooth them first.
Here we adopt the smoothing technique discussed in section \ref{sec:nonsmooth}; we approximate $f_i(x)$ by the smooth function
\bea\label{smooth-fi}f_i^\rho(x):=\max\{\langle x-c^i, y\rangle -\frac{\rho}{2}\|y\|^2 : \|y\|\leq 1\},\eea where $\rho>0$
is a smoothness parameter. The gradient of $f_i^\rho$, $\nabla f_i^\rho(x)= y_i^*,$ where $y_i^*$ is the optimal
solution to the optimization problem in \eqref{smooth-fi}. It is easy to show that $y_i^*=\frac{x-c^i}{\max\{\rho,\|x-c^i\|\}}.$ Moreover, $\nabla f_i^\rho(x)$ is Lipschitz continuous with constant $L(f_i^\rho)=1/\rho$.
Now we can apply MSA, FaMSA and FaMSA-s to solve \bea\label{prob:weber-smooth}\min\sum_{i=1}^K f_i^\rho(x).\eea
The $i$-th subproblem in all of these algorithms corresponds to solving the following problem:
\bea\label{subproblem-weber}p_i(w_{(k)}^i,\ldots,w_{(k)}^i):=\arg\min_u f_i^\rho(u)+\sum_{j=1,j\neq i}^K
\left(f_j^\rho(w_{(k)}^i)+\langle\nabla f_j^\rho(w_{(k)}^i),u-w_{(k)}^i\rangle+\oneotwomu\|u-w_{(k)}^i\|^2\right).\eea
It is easy to check that the optimal solution to problem \eqref{subproblem-weber} is given by \beaa u^*:=\left\{\ba{ll}
\displaystyle c^i+\frac{\rho(K-1)}{\mu+\rho(K-1)}(z_{(k)}^i-c^i), & \mbox{if } \|z_{(k)}^i-c^i\|\leq\rho+\frac{\mu}{K-1} \\
\displaystyle c^i+\frac{(K-1)\|z_{(k)}^i-c^i\|-\mu}{(K-1)\|z_{(k)}^i-c^i\|}(z_{(k)}^i-c^i), & \mbox{if } \|z_{(k)}^i-c^i\|>\rho+\frac{\mu}{K-1},  \ea\right.\eeaa
where $$z_{(k)}^i:=w_{(k)}^i-\frac{\mu}{K-1}\sum_{j=1,j\neq i}^K \frac{w_{(k)}^i-c^j}{\max\{\rho,\|w_{(k)}^i-c^j\|\}}.$$

If we choose the doubly stochastic matrix $D^{(k)}$ to be $D^{(k)}:=1/Kee^\top$ in MSA as we do in FaMSA-s, all $w_{(k)}^i$'s are the same in MSA as they are in FaMSA-s. Hence, computing $x_{(k)}^i$, for $i=1,\ldots,K$ in both algorithms can be done efficiently as follows.
\bea\label{alg:num-imple-FW} \left\{\ba{lll} \hat{z}_{(k)} & = & \sum_{j=1}^K \frac{w_{(k)}-c^j}{\max\{\rho,\|w_{(k)}-c^j\|\}} \\ z_{(k)}^i & = & w_{(k)} - \frac{\mu}{K-1}(\hat{z}-\frac{w_{(k)}-c^i}{\max\{\rho,\|w_{(k)}-c^i\|\}}), \forall i = 1,\ldots,K \\ x_{(k)}^i & = & c^i + (1-\frac{\mu}{\max\{(K-1)\|z_{(k)}^i-c^i\|,\mu+\rho(K-1)\}})(z_{(k)}^i-c^i), \forall i = 1,\ldots,K \ea\right.\eea

We compared the performance of MSA and FaMSA-s with the classical gradient method (Grad) and Nesterov's accelerated gradient method (Nest) for solving \eqref{prob:weber-smooth}. The classical gradient method for solving \eqref{prob:weber-smooth} with step size $\tau>0$ is: \beaa x^{k+1} = x^k - \tau \sum_{j=1}^K \nabla f_j^\rho(x^k).  \eeaa The variant of Nesterov's accelerated gradient method that we used is the following:
\beaa \left\{ \ba{lll} x^{k} & = & y^{k-1} - \tau \sum_{j=1}^K \nabla f_j^\rho(y^{k-1}) \\ y^{k} & = & x^{k} - \frac{k-1}{k+2}(x^k -x^{k-1}).  \ea \right. \eeaa

We created random problems to test the performance of MSA, FaMSA-s, Grad and Nest as follows. Vectors $c^i\in\br^n,i=1,\ldots,K$ were created
with i.i.d. Gaussian entries from $\mathcal{N}(0,n)$. The seed for generating random numbers in MATLAB was set to 0. We set the smoothness parameter $\rho$ equal to $10^{-3}$.
The initial points $x^i, i=1,\ldots,K$ were set to the average of all of the $c^i$'s, i.e., $x_{(0)}^i=\frac{1}{K}\sum_{i=1}^K c^i.$
We chose $D^{(k)}_{ij}=1/K,i,j=1,\ldots,K$ for all $k$ in MSA. To compare the number of iterations needed by MSA and FaMSA-s, we first solved \eqref{prob:weber} by Mosek \cite{Mosek} after converting it into an
SOCP problem to get the optimal solution $x^*$, and then terminated MSA, FaMSA-s, Grad and Nest when the relative error of the objective function
value at the $k$-th iterate,
\beaa relerr:=\frac{|\min_{i=1,\ldots,K}F(x_{(k)}^i)-F(x^*)|}{F(x^*)},\eeaa was less than $10^{-6}$.
We tested the performance of these four solvers for different choices of $\tau$, which is the step size for Grad
and Nest. Note that since the $w^i$'s are the same in MSA with $D^{(k)}=\frac{1}{K}ee^\top$ for all $k$ and in FaMSA-s, these two methods can be viewed as linearization methods in which the single function $\sum_{j=1,j\neq i}^K f_j(x)$ is linearized at the point $w$ with only one proximal term $\frac{K-1}{2\mu}\|x-w\|$ in the $i$-th subproblem. So the step size for MSA and FaMSA-s is $\mu/(K-1)$. Hence, the parameter $\mu$ for MSA and FaMSA-s was set to $\mu=\tau(K-1)$ in our numerical tests.

Our results are presented in Table \ref{tab:weber}. The CPU times reported are in seconds. These results show that for the F-W problem, our implementations of MSA and FaMSA-s take roughly between two and three times as much time to solve each problem as taken by Grad and Nest, respectively. This is not surprising since it is clear that the computation of each set of $K$ vectors $z_{(k)}^i$ and $x_{(k)}^i$ for $i=1,\ldots,K$ in \eqref{alg:num-imple-FW} is roughly comparable to a single computation of the gradient, i.e., the $K$ gradients of $f_i^\rho(x)$, for $i=1,\ldots,K$. Moreover, for the simple F-W objective function, not much is gained by minimizing only one out of the $K$ individual functions $f_i^\rho(x)$, $i=1,\ldots,K$, when $K$ is large as it is in our tests. Note that the number of iterations required by MSA and Grad were exactly the same on our set of test problems. When $K$ is of a moderate size and the individual functions are more complicated, MSA should require fewer iterations than Grad.

\begin{table}[ht]{\small
\begin{center}\caption{Comparison of MSA, FaMSA-v, Grad and Nest on solving Fermat-Weber problem \eqref{prob:weber-smooth}}\label{tab:weber}
\begin{tabular}{|c c |c | c c c | c c c | c c c | c c c |}\hline
\multicolumn{2}{|c|}{Problem} & \multicolumn{1}{|c|}{Mosek} & \multicolumn{3}{c|}{MSA} &
\multicolumn{3}{c|}{FaMSA-s} & \multicolumn{3}{c|}{Grad} & \multicolumn{3}{c|}{Nest} \\\hline

$n$ & $K$ & time & iter & relerr & time & iter & relerr & time & iter & relerr & time & iter & relerr & time \\\hline

\multicolumn{15}{|c|}{$\tau = 0.001$} \\\hline
50 & 50 & 0.85 & 500 & 4.1e-05 & 0.73 & 107 & 8.4e-07 & 0.17 & 500 & 4.1e-05 & 0.21 & 109 & 9.0e-07 & 0.05 \\\hline
50 & 100 & 3.40 & 500 & 5.6e-06 & 1.44 & 69 & 9.9e-07 & 0.21 & 500 & 5.6e-06 & 0.42 & 72 & 8.5e-07 & 0.07  \\\hline
50 & 200 & 0.96 & 427 & 9.9e-07 & 2.44 & 47 & 8.8e-07 & 0.28 & 427 & 9.9e-07 & 0.69 & 49 & 8.6e-07 & 0.09  \\\hline
100 & 100 & 1.78 & 500 & 8.9e-06 & 1.68 & 94 & 9.8e-07 & 0.33 & 500 & 8.9e-06 & 0.51 & 97 & 9.1e-07 & 0.10 \\\hline
100 & 200 & 4.48 & 500 & 1.6e-06 & 3.35 & 60 & 9.2e-07 & 0.42 & 500 & 1.6e-06 & 1.00 & 62 & 9.2e-07 & 0.13 \\\hline
100 & 400 & 9.40 & 198 & 1.0e-06 & 2.68 & 34 & 9.5e-07 & 0.47 & 198 & 1.0e-06 & 0.79 & 36 & 9.2e-07 & 0.15 \\\hline
200 & 200 & 22.22 & 500 & 2.3e-06 & 4.36 & 75 & 9.9e-07 & 0.67 & 500 & 2.3e-06 & 1.39 & 77 & 9.9e-07 & 0.22 \\\hline
200 & 400 & 45.55 & 275 & 1.0e-06 & 4.81 & 41 & 9.9e-07 & 0.73 & 275 & 1.0e-06 & 1.54 & 43 & 9.8e-07 & 0.25 \\\hline
200 & 800 & 100.15 & 41 & 1.0e-06 & 1.44 & 15 & 9.7e-07 & 0.54 & 41 & 1.0e-06 & 0.46 & 16 & 9.8e-07 & 0.18 \\\hline
300 & 300 & 102.64 & 419 & 1.0e-06 & 6.73 & 52 & 9.9e-07 & 0.85 & 419 & 1.0e-06 & 2.22 & 54 & 9.9e-07 & 0.29 \\\hline
300 & 600 & 194.99 & 24 & 1.0e-06 & 0.79 & 11 & 9.9e-07 & 0.37 & 24 & 1.0e-06 & 0.26 & 12 & 9.9e-07 & 0.14 \\\hline
300 & 1200 & 401.54 & 1 & 5.8e-07 & 0.08 & 1 & 5.8e-07 & 0.08 & 1 & 5.8e-07 & 0.03 & 1 & 5.8e-07 & 0.03 \\\hline
\multicolumn{15}{|c|}{$\tau = 0.01$} \\\hline
50 & 50 & 0.84 & 238 & 9.9e-07 & 0.36 & 32 & 8.3e-07 & 0.06 & 238 & 9.9e-07 & 0.11 & 34 & 7.5e-07 & 0.02 \\\hline
50 & 100 & 3.36 & 93 & 9.9e-07 & 0.29 & 20 & 9.6e-07 & 0.07 & 93 & 9.9e-07 & 0.08 & 22 & 7.6e-07 & 0.03 \\\hline
50 & 200 & 0.96 & 42 & 9.9e-07 & 0.26 & 13 & 9.2e-07 & 0.09 & 42 & 9.9e-07 & 0.07 & 15 & 5.9e-07 & 0.03 \\\hline
100 & 100 & 1.78 & 160 & 1.0e-06 & 0.55 & 28 & 9.0e-07 & 0.11 & 160 & 1.0e-06 & 0.17 & 30 & 8.1e-07 & 0.04 \\\hline
100 & 200 & 4.48 & 62 & 9.8e-07 & 0.43 & 17 & 9.2e-07 & 0.13 & 62 & 9.8e-07 & 0.13 & 19 & 7.5e-07 & 0.05 \\\hline
100 & 400 & 9.46 & 20 & 9.5e-07 & 0.28 & 9 & 9.1e-07 & 0.13 & 20 & 9.5e-07 & 0.09 & 10 & 9.2e-07 & 0.05 \\\hline
200 & 200 & 22.37 & 91 & 1.0e-06 & 0.81 & 22 & 9.2e-07 & 0.21 & 91 & 1.0e-06 & 0.26 & 23 & 1.0e-06 & 0.07 \\\hline
200 & 400 & 45.56 & 28 & 9.7e-07 & 0.50 & 11 & 9.9e-07 & 0.21 & 28 & 9.7e-07 & 0.16 & 13 & 8.4e-07 & 0.08 \\\hline
200 & 800 & 99.38 & 4 & 1.0e-06 & 0.16 & 4 & 8.6e-07 & 0.16 & 4 & 1.0e-06 & 0.05 & 4 & 9.4e-07 & 0.05 \\\hline
300 & 300 & 100.48 & 42 & 9.9e-07 & 0.69 & 15 & 9.3e-07 & 0.26 & 42 & 9.9e-07 & 0.23 & 16 & 9.5e-07 & 0.09 \\\hline
300 & 600 & 194.88 & 3 & 9.7e-07 & 0.11 & 3 & 9.4e-07 & 0.11 & 3 & 9.7e-07 & 0.04 & 3 & 9.6e-07 & 0.04 \\\hline
300 & 1200 & 402.16 & 1 & 5.4e-07 & 0.08 & 1 & 5.4e-07 & 0.08 & 1 & 5.4e-07 & 0.03 & 1 & 5.4e-07 & 0.03 \\\hline
\multicolumn{15}{|c|}{$\tau = 0.1$} \\\hline
50 & 50 & 0.84 & 23 & 9.5e-07 & 0.05 & 9 & 3.4e-07 & 0.03 & 23 & 9.4e-07 & 0.02 & 10 & 5.4e-07 & 0.01 \\\hline
50 & 100 & 3.41 & 9 & 7.7e-07 & 0.04 & 5 & 6.1e-07 & 0.03 & 9 & 7.7e-07 & 0.02 & 6 & 5.1e-07 & 0.01  \\\hline
50 & 200 & 0.95 & 4 & 5.3e-07 & 0.04 & 3 & 2.9e-07 & 0.03 & 4 & 5.2e-07 & 0.01 & 3 & 1.0e-06 & 0.01 \\\hline
100 & 100 & 1.80 & 16 & 8.6e-07 & 0.07 & 8 & 3.6e-07 & 0.04 & 16 & 8.6e-07 & 0.02 & 9 & 4.1e-07 & 0.02 \\\hline
100 & 200 & 4.48 & 6 & 8.3e-07 & 0.05 & 4 & 7.2e-07 & 0.04 & 6 & 8.3e-07 & 0.02 & 5 & 5.2e-07 & 0.02 \\\hline
100 & 400 & 9.40 & 2 & 6.4e-07 & 0.04 & 2 & 4.2e-07 & 0.04 & 2 & 6.4e-07 & 0.02 & 2 & 6.4e-07 & 0.02 \\\hline
200 & 200 & 22.25 & 9 & 9.4e-07 & 0.09 & 6 & 5.8e-07 & 0.07 & 9 & 9.3e-07 & 0.03 & 6 & 9.4e-07 & 0.02 \\\hline
200 & 400 & 45.61 & 3 & 7.9e-07 & 0.07 & 3 & 5.0e-07 & 0.07 & 3 & 7.9e-07 & 0.02 & 3 & 6.9e-07 & 0.03 \\\hline
200 & 800 & 99.77 & 1 & 5.0e-07 & 0.05 & 1 & 5.0e-07 & 0.05 & 1 & 5.0e-07 & 0.02 & 1 & 5.0e-07 & 0.02 \\\hline
300 & 300 & 100.37 & 4 & 9.9e-07 & 0.08 & 4 & 6.7e-07 & 0.08 & 4 & 9.9e-07 & 0.03 & 4 & 8.4e-07 & 0.03 \\\hline
300 & 600 & 197.72 & 1 & 7.0e-07 & 0.05 & 1 & 7.0e-07 & 0.05 & 1 & 7.0e-07 & 0.02 & 1 & 7.0e-07 & 0.02 \\\hline
300 & 1200 & 412.49 & 1 & 2.1e-07 & 0.08 & 1 & 2.1e-07 & 0.08 & 1 & 2.1e-07 & 0.03 & 1 & 2.1e-07 & 0.03 \\\hline
\end{tabular}
\end{center}}
\end{table}

The purpose of this set of tests was not to demonstrate any advantage that our algorithms might have over gradient methods. Rather, they were performed to validate our algorithms and show that the accelerated variants like algorithm Nest can reduce the number of iterations required to solve problems of the form \eqref{prob:convex-K-sum-original}. This is quite clear from the results reported in Table \ref{tab:weber}. We further note that FaMSA-s often takes one to three fewer iterations than Nest. Note that for some problems, the multiple splitting algorithm took only one iteration to converge. The reason was that for these problems, the number of points was much larger than the dimension of the space. Therefore, the points were very compact and fairly uniformly distributed around the initial point; hence that point was quite likely to be very close to the optimal solution.

\subsection{An image deblurring problem}

In this section, we report the results of applying our multiple splitting algorithms to a benchmark total variation and wavelet-based image deblurring problem from
\cite{Figueiredo-Nowak-03}.
In this problem, the original image is the well-known Cameraman image of size $256\times 256$ and the observed image is obtained after imposing a uniform blur of
size $9\times 9$ (denoted by the operator $A$) and Gaussian noise (generated by the function {\it randn}
in MATLAB with a seed of 0 and a standard deviation of $0.56$).
Since the vector of coefficients of the wavelet transform of the image is sparse in this problem and the total variation norm of the image is expected
to be small, one can try to reconstruct the image $x$ from the observed image $b$ by solving the problem:
\bea\label{prob:TV-wavelet-deblur} \min \quad \alpha \TV(x) + \beta \|\Phi x\|_1 + \half\|Ax-b\|_2^2, \eea
where $\TV(x):=\sum_{ij}\sqrt{(x_{i+1,j}-x_{ij})^2+(x_{ij}-x_{i,j+1})^2}$ is the total variation of $x$, $\Phi$ is the wavelet transform, $A$ denotes
the deblurring kernel and $\alpha>0$, $\beta>0$ are weighting parameters. Problem \eqref{prob:TV-wavelet-deblur} involves minimizing the sum of three
convex functions with $f_1(x)=\alpha \TV(x)$, $f_2(x)=\beta \|\Phi x\|_1$ and $f_3(x)=\half\|Ax-b\|_2^2$.

To apply our multiple splitting algorithms to solve
\eqref{prob:TV-wavelet-deblur}, our theory
requires all the functions to be smooth functions. So we needed to smooth the $\TV$ and the $\ell_1$ functions first.
We adopted the following way to smooth the $\TV$ function, widely used in the literature for doing this:
\beaa f_1^\delta(x) := \alpha\sum_{ij}\sqrt{(x_{i+1,j}-x_{ij})^2+(x_{ij}-x_{i,j+1})^2+\delta}. \eeaa The $\ell_1$ function was smoothed in the way described in
Section \ref{sec:nonsmooth}: \beaa f_2^\sigma(x) := \beta\max_u\{\langle \Phi x,u\rangle - \frac{\sigma}{2}\|u\|^2 : \|u\|_\infty\leq 1\}. \eeaa Thus,
the smooth version of problem \eqref{prob:TV-wavelet-deblur} was:
\bea \label{prob:TV-wavelet-deblur-smooth} \min_x \quad f_1^\delta(x) + f_2^\sigma(x) + f_3(x). \eea

However, when we applied our multiple splitting algorithms
to \eqref{prob:TV-wavelet-deblur-smooth}, we actually performed the following computation on the $k$-th iteration:
\bea\left\{\label{alg:msa-tv-wavelet-deblur}\ba{lll}x^{k+1} & := & \arg\min_x f_1(x) + \langle\nabla f_2^\sigma(w^k), x-w^k\rangle+\frac{1}{2\mu}\|x-w^k\|^2+
\langle\nabla f_3(w^k),x-w^k\rangle+\frac{1}{2\mu}\|x-w^k\|^2 \\
                                                    y^{k+1} & := & \arg\min_y f_2^\sigma(y) + \langle\nabla f_1^\delta(w^k),y-w^k\rangle+\frac{1}{2\mu}\|y-w^k\|^2+
                                                    \langle\nabla f_3(w^k),y-w^k\rangle+\frac{1}{2\mu}\|y-w^k\|^2\\
                                                    z^{k+1} & := & \arg\min_z f_3(z) + \langle\nabla f_1^\delta(w^k),z-w^k\rangle+\frac{1}{2\mu}\|z-w^k\|^2+
                                                    \langle\nabla f_2^\sigma(w^k),z-w^k\rangle+\frac{1}{2\mu}\|z-w^k\|^2 \\
                                                    w^{k+1} & := & (x^{k+1}+y^{k+1}+z^{k+1})/3. \ea\right.\eea
Note that in \eqref{alg:msa-tv-wavelet-deblur}, when we linearized the $\TV$ function, we used the smoothed $\TV$ function $f_1^\delta(\cdot)$,
i.e., we computed the
gradient of $f_1^\delta(\cdot)$. But when we solved the first subproblem, we used the nonsmooth $\TV$ function $f_1(\cdot)$, because there are efficient algorithms for solving this nonsmooth problem. Specifically, this subproblem can be reduced to: \beaa x^{k+1}:=\arg\min_x \frac{\alpha\mu}{2}\TV(x) + \half\|x-
\left(w^k-\mu(\nabla f_2^\sigma(w^k)+\nabla f_3(w^k)\right)/2)\|^2, \eeaa which
is a standard TV-denoising problem. In our tests, we perform 10 iterations of the algorithm proposed by Chambolle in \cite{Chambolle-04} to approximately solve this problem.
The second subproblem in \eqref{alg:msa-tv-wavelet-deblur} can be reduced to:
\bea\label{prob:Wavelet-prox}y^{k+1}:= \arg\min_y \frac{\mu}{2}f_2^\sigma(y) + \half\|y-(w^k-\mu\left(\nabla f_1^\delta(w^k)+\nabla f_3(w^k)\right)/2)\|^2. \eea
It is easy to check that the solution of \eqref{prob:Wavelet-prox} is given by:
\beaa y^{k+1}:=\Phi^\top\left(\Phi \bar{w}^k - \frac{\mu\beta}{2}\tilde{w}^k\right) \eeaa
where $(\tilde{w}^k)_j = \max\{-1,\min\{1,\frac{2(\Phi\bar{w}^k)_j}{2\sigma+\beta\mu}\}\}$ and
$\bar{w}^k=w^k-\mu\left(\nabla f_1^\delta(w^k)+\nabla f_3(w^k)\right)/2$. The third subproblem in \eqref{alg:msa-tv-wavelet-deblur} corresponds to solving
the following linear system: \beaa (A^\top A+2/\mu I)z = A^\top b - \nabla f_1^\delta(w^k) + 2/\mu w^k - \nabla f_2^\delta(w^k).  \eeaa
Solving this linear system is easy since the operator $A$ has a special structure and thus $(A^\top A+ 2/\mu I)$ can be inverted efficiently.

\begin{figure}
\centering \subfigure{\includegraphics[scale=0.5]{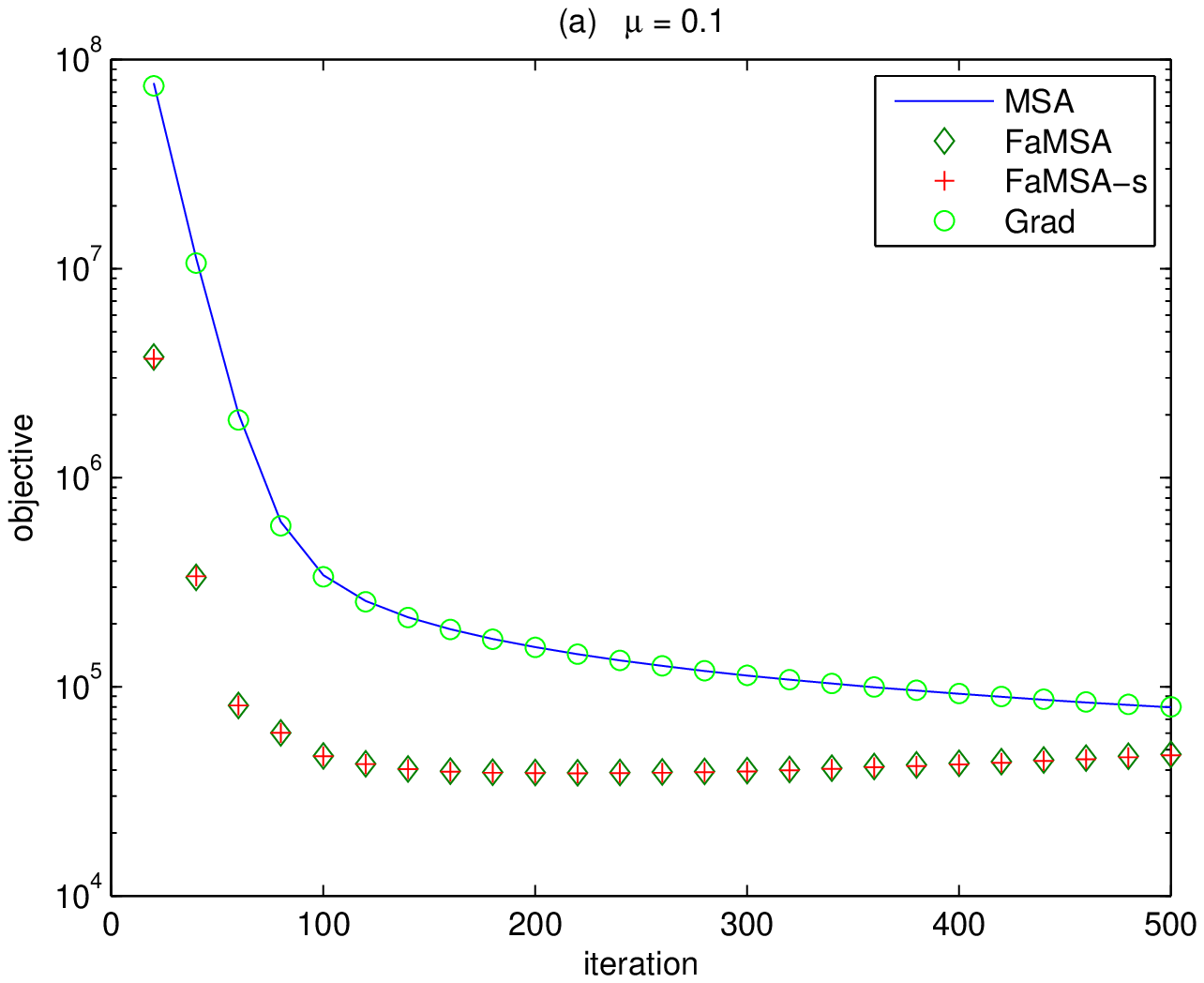}}\hspace{-0.5cm}
\centering \subfigure{\includegraphics[scale=0.5]{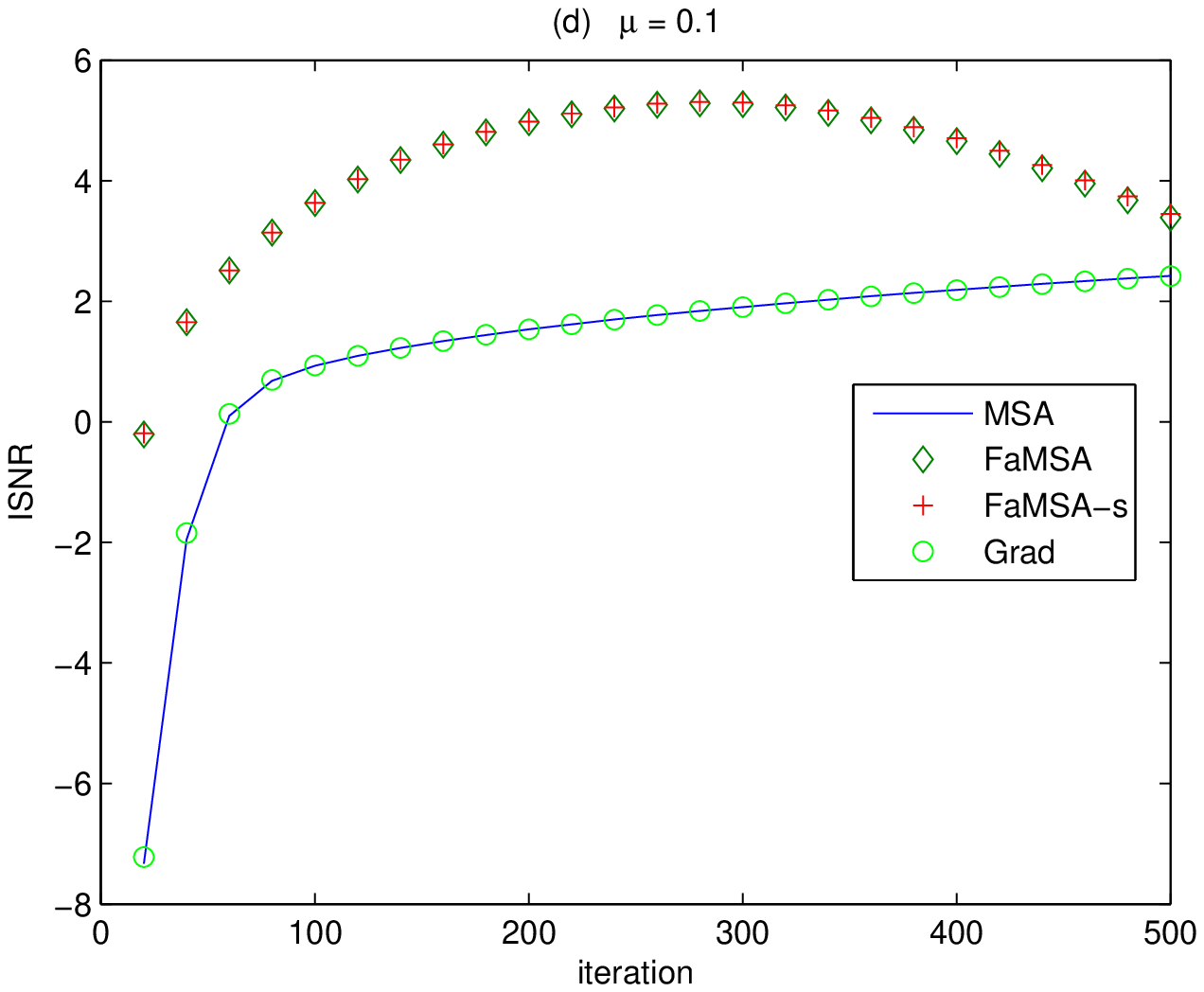}}\\
\centering \subfigure{\includegraphics[scale=0.5]{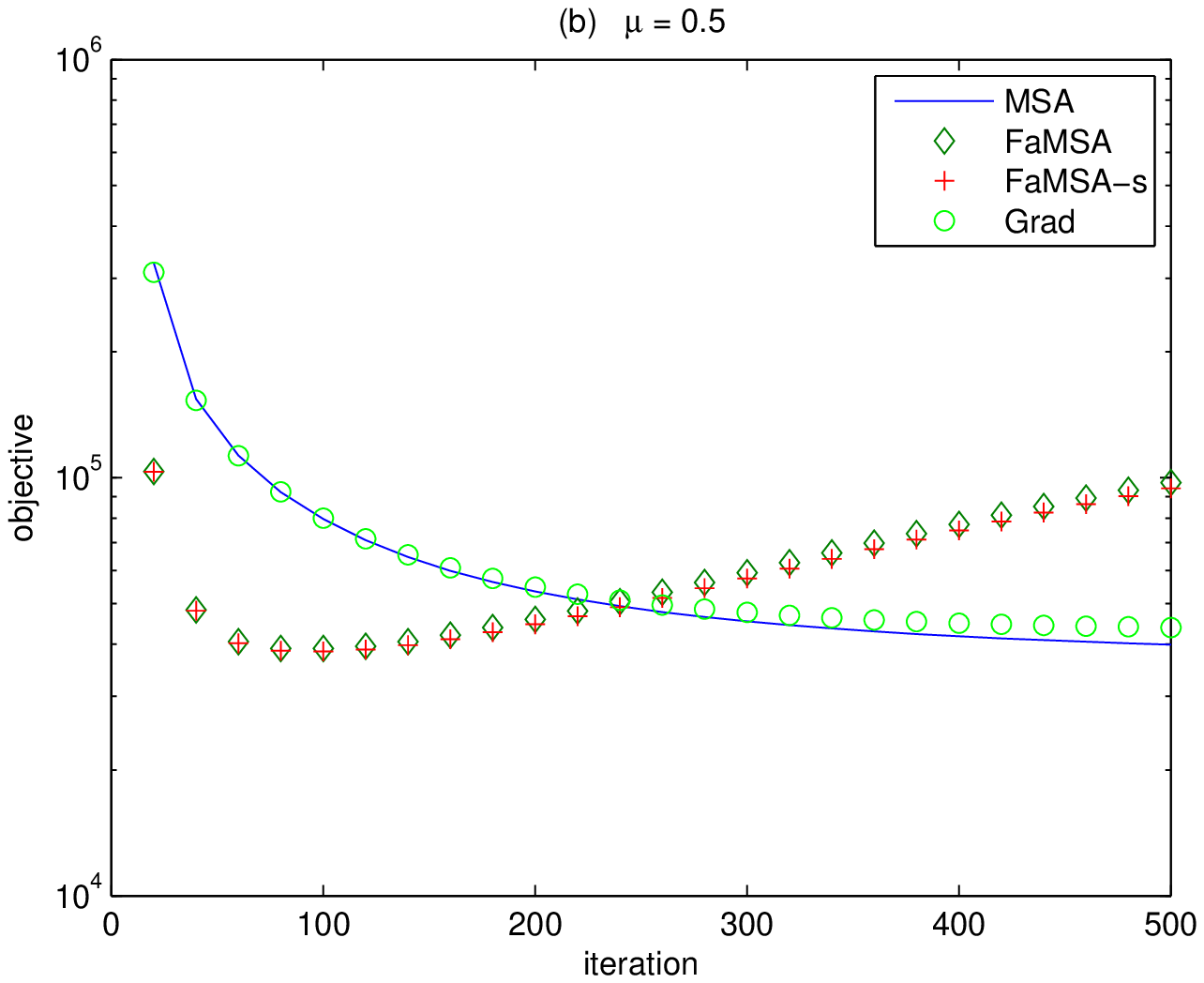}} \hspace{-0.5cm}
\centering \subfigure{\includegraphics[scale=0.5]{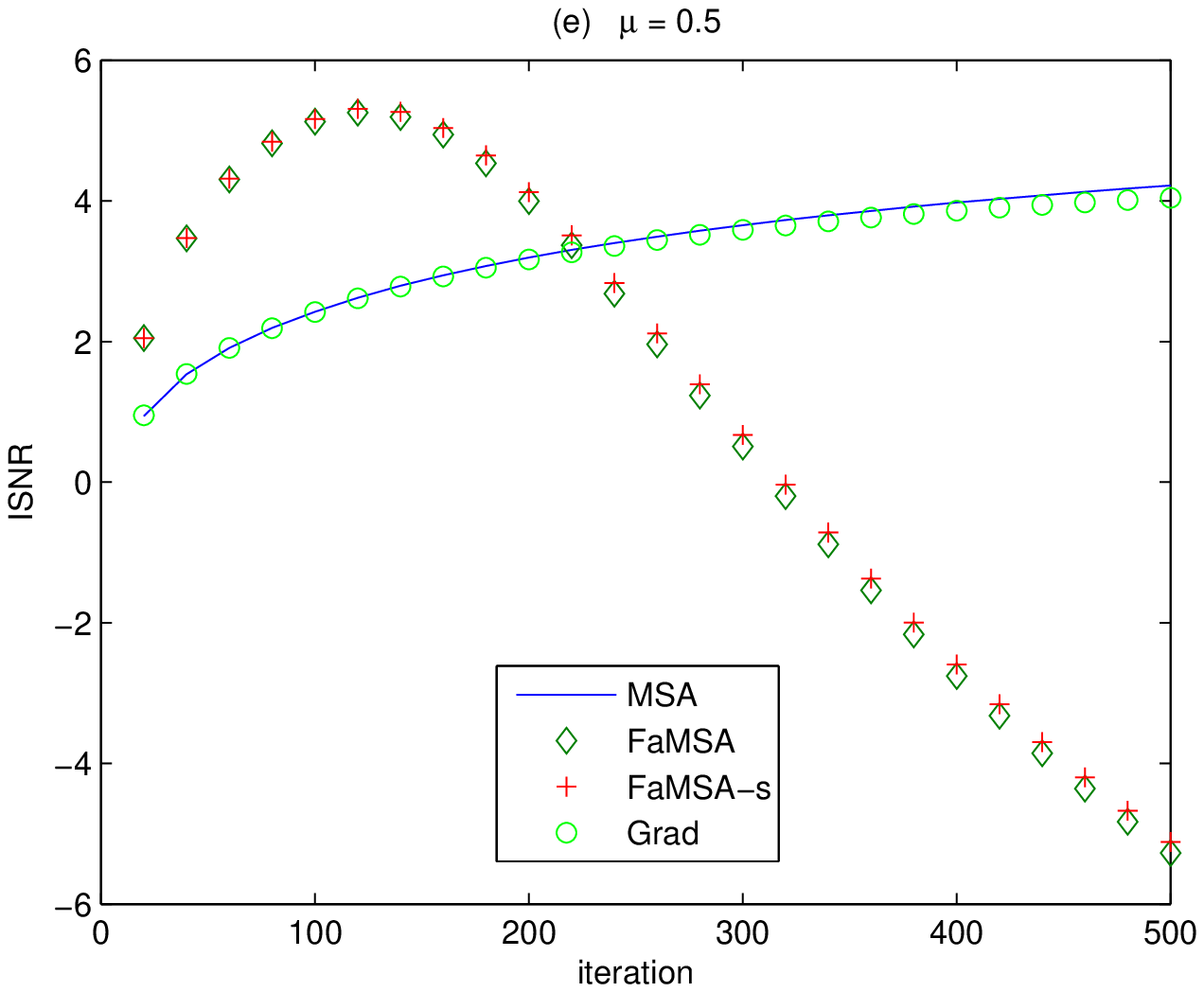}}\\
\centering \subfigure{\includegraphics[scale=0.5]{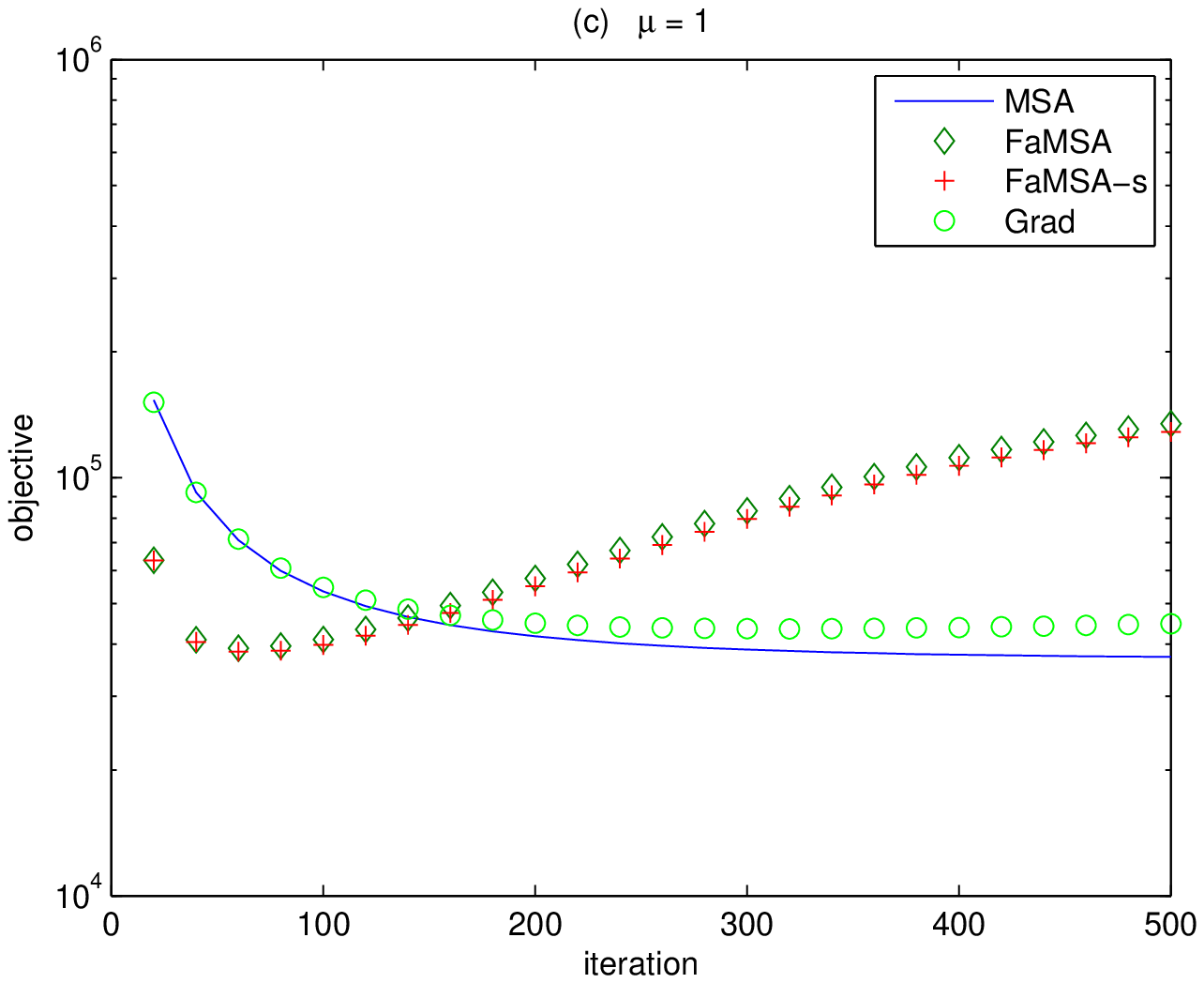}} \hspace{-0.5cm}
\centering \subfigure{\includegraphics[scale=0.5]{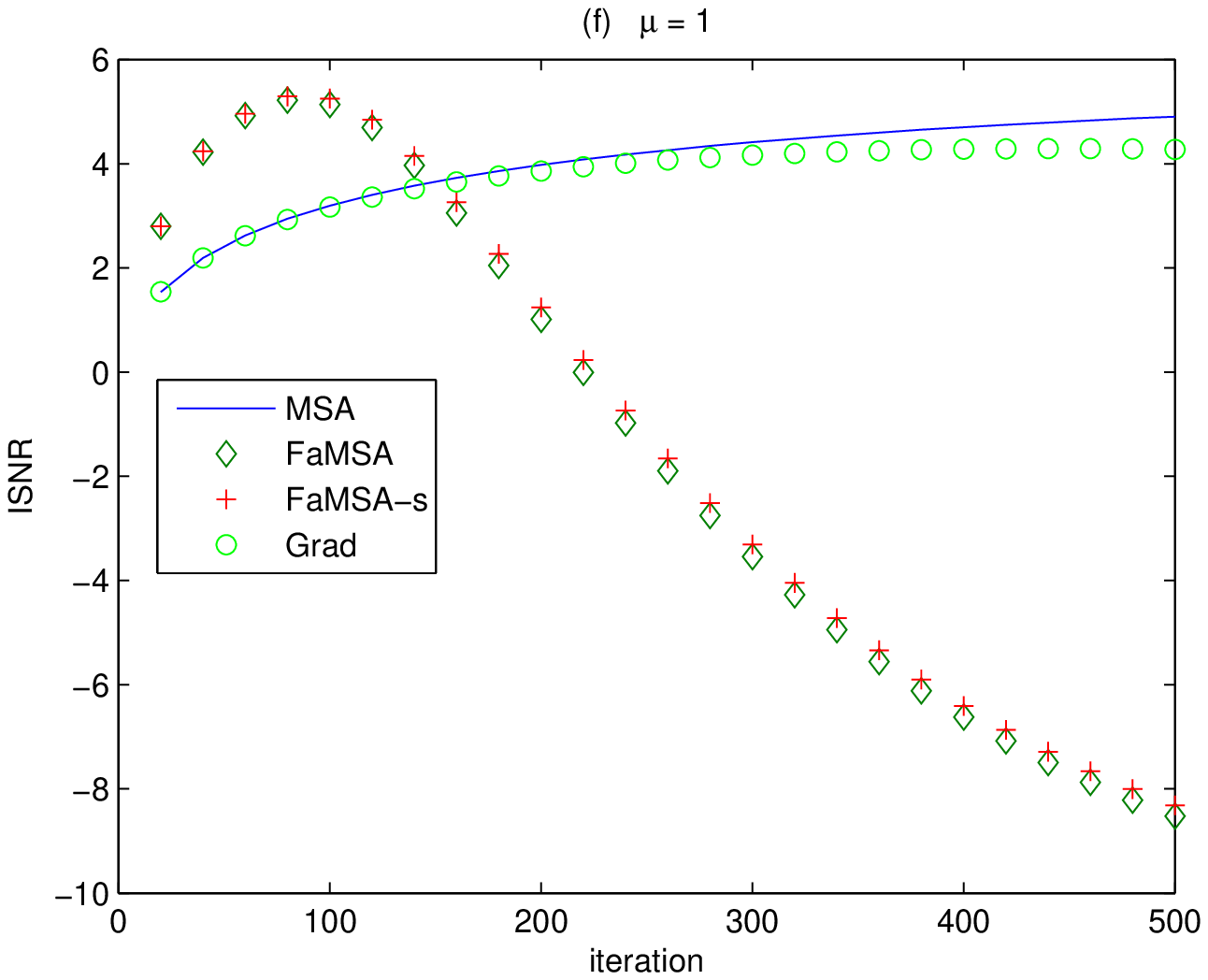}}
\caption{Comparison of MSA, FaMSA-s and Grad for different $\mu$.}\label{fig:compare-deblur}
\end{figure}

In our tests, we set $\alpha=0.001$, $\beta=0.035$ and used smoothing parameters $\delta=\sigma=10^{-4}$. The initial points were all set equal to $0$.
We compared the performance of MSA, FaMSA, FaMSA-s and Grad for different $\mu$ and step sizes $\tau$. In these comparisons, we simply terminated the codes after 500 iterations. The objective function value and the improvement signal noise ratio (ISNR) at different iterations are reported in Table \ref{tab:TV-deblur-test}. The ISNR is defined as $ISNR:= 10\log_{10}\frac{\|x-b\|^2}{\|x-\bar{x}\|^2}$, where $x$ is the reconstructed image and $\bar{x}$ is the true image. As we did for F-W problem, we always used $\mu=\tau(K-1)$ and since there were three functions in this problem, we used $\mu=2\tau$. For large $\mu$, we did not report the results for all of the iterations since the comparisons are quite clear from the selected iterations. See Figure \ref{fig:compare-deblur} for additional and more complete comparisons. We make the following observations from Table \ref{tab:TV-deblur-test}.
For $\mu=0.1$, FaMSA-s achieved the best objective function value in about 200 iterations and 152 CPU seconds. The best ISNR was also achieved by FaMSA-s, in about 300 iterations and 227 seconds. MSA and Grad were not able to obtain an acceptable solution in 500 iterations. In fact, they were only able to reduce the objective function to twice the near-optimal value of $3.86\times 10^4$ achieved by FaMSA-s. For $\mu=0.5$, FaMSA-s achieved the best objective function value and ISNR in 100 iterations and 76 seconds and 125 iterations and 94 seconds, respectively. Again, MSA and Grad did not achieve acceptable results even after 500 iterations. For $\mu=1$, MSA achieved the best objective function value, $3.73\times 10^4$, after 500 iterations and 349 CPU seconds, while the best ISNR was achieved by FaMSA-s in 80 iterations and 61 seconds. Also, the best objective function value achieved by FaMSA-s was at the 60-th iteration after only 47 CPU seconds. We also note that for $\mu=0.1,0.5$ and $1$, MSA was always better than Grad and FaMSA-s was always slightly better than FaMSA. Another observation was that MSA always decreased the objective function value for $\mu=0.1,0.5$ and $1$, while FaMSA and FaMSA-s always achieved near-optimal results in a relatively small number of iterations and then started getting worse. However, in practice, one would always terminate FaMSA and FaMSA-s once the objective function value started increasing. For $\mu=5$, MSA gave very good results while the other three solvers diverged immediately. Specifically, the best objective function value $3.73\times 10^4$ was achieved by MSA in 120 iterations and 80 CPU seconds, and the best ISNR was achieved by MSA in 200 iterations and 132 CPU seconds. Thus, based on these observations, we conclude that FaMSA-s attains a nearly optimal solution very quickly for small $\mu$ while MSA is more stable for large $\mu$.


\begin{table}[ht]{\small
\begin{center}\caption{Comparison of MSA, FaMSA, FaMSA-s and Grad on solving TV-deblurring problem}\label{tab:TV-deblur-test}
\begin{tabular}{|c|cc|cc|cc|cc|}\hline
& \multicolumn{2}{|c|}{MSA} & \multicolumn{2}{c|}{FaMSA} & \multicolumn{2}{c|}{FaMSA-s} & \multicolumn{2}{|c|}{Grad} \\\hline

Iter & obj & ISNR & obj & ISNR & obj & ISNR & obj & ISNR \\\hline

\multicolumn{9}{|c|}{$\mu = 0.1, \tau = 0.05$ } \\\hline

100 & 3.42e+005 &     0.9311 & 4.67e+004 &     3.6310 & 4.66e+004 &     3.6332 & 3.36e+005 &     0.9344 \\\hline

200 & 1.55e+005 &     1.5340 & 3.89e+004 &     4.9693 & {\bf 3.86e+004} &     4.9821 & 1.55e+005 &     1.5341 \\\hline

300 & 1.13e+005 &     1.9057 & 3.98e+004 &     5.2695 & 3.94e+004 &     {\bf 5.2989} & 1.13e+005 &     1.9043 \\\hline

400 & 9.25e+004 &     2.1905 & 4.30e+004 &     4.6587 & 4.26e+004 &     4.7075 & 9.28e+004 &     2.1871 \\\hline

500 & 7.97e+004 &     2.4235 & 4.76e+004 &     3.3881 & 4.70e+004 &     3.4500 & 8.02e+004 &     2.4175 \\\hline

\multicolumn{9}{|c|}{$\mu = 0.5, \tau = 0.25$ } \\\hline

25 & 2.41e+005 &     1.1343 & 7.70e+004 &     2.4777 & 7.69e+004 &     2.4784 & 2.36e+005 &     1.1408 \\\hline

50 & 1.29e+005 &     1.7359 & 4.31e+004 &     3.9343 & 4.28e+004 &     3.9416 & 1.29e+005 &     1.7376 \\\hline

75 & 9.66e+004 &     2.1260 & 3.92e+004 &     4.7122 & 3.88e+004 &     4.7324 & 9.67e+004 &     2.1250 \\\hline

100 & 7.96e+004 &     2.4243 & 3.90e+004 &     5.1257 & {\bf 3.84e+004} &     5.1638 & 8.00e+004 &     2.4198 \\\hline

125 & 6.92e+004 &     2.6659 & 3.97e+004 &     5.2558 & 3.90e+004 &     {\bf 5.3160} & 6.98e+004 &     2.6569 \\\hline

150 & 6.21e+004 &     2.8682 & 4.12e+004 &     5.0880 & 4.04e+004 &     5.1737 & 6.30e+004 &     2.8538 \\\hline

175 & 5.71e+004 &     3.0416 & 4.33e+004 &     4.6478 & 4.23e+004 &     4.7576 & 5.82e+004 &     3.0207 \\\hline

200 & 5.34e+004 &     3.1928 & 4.58e+004 &     3.9964 & 4.46e+004 &     4.1258 & 5.48e+004 &     3.1646 \\\hline

225 & 5.06e+004 &     3.3266 & 4.86e+004 &     3.2006 & 4.73e+004 &     3.3442 & 5.22e+004 &     3.2902 \\\hline

250 & 4.85e+004 &     3.4463 & 5.18e+004 &     2.3223 & 5.03e+004 &     2.4758 & 5.03e+004 &     3.4009 \\\hline

275 & 4.67e+004 &     3.5545 & 5.54e+004 &     1.4132 & 5.37e+004 &     1.5723 & 4.88e+004 &     3.4991 \\\hline

300 & 4.54e+004 &     3.6529 & 5.93e+004 &     0.5078 & 5.74e+004 &     0.6705 & 4.76e+004 &     3.5869 \\\hline

500 & 3.99e+004 &     4.2186 & 9.74e+004 &    -5.2730 & 9.43e+004 &    -5.1193 & 4.38e+004 &     4.0416 \\\hline

\multicolumn{9}{|c|}{$\mu = 1, \tau = 0.5$ } \\\hline

20 & 1.53e+005 &     1.5382 & 6.35e+004 &     2.7991 & 6.33e+004 &     2.8006 & 1.51e+005 &     1.5443 \\\hline

40 & 9.23e+004 &     2.1927 & 4.10e+004 &     4.2214 & 4.05e+004 &     4.2361 & 9.22e+004 &     2.1932 \\\hline

60 & 7.09e+004 &     2.6220 & 3.91e+004 &     4.9205 & {\bf 3.84e+004} &     4.9591 & 7.13e+004 &     2.6158 \\\hline

80 & 5.99e+004 &     2.9413 & 3.96e+004 &     5.2175 & 3.86e+004 &     {\bf 5.2890} & 6.08e+004 &     2.9258 \\\hline

100 & 5.34e+004 &     3.1933 & 4.10e+004 &     5.1371 & 3.98e+004 &     5.2488 & 5.47e+004 &     3.1664 \\\hline

120 & 4.93e+004 &     3.4003 & 4.33e+004 &     4.6922 & 4.19e+004 &     4.8439 & 5.10e+004 &     3.3597 \\\hline

140 & 4.64e+004 &     3.5751 & 4.62e+004 &     3.9649 & 4.45e+004 &     4.1489 & 4.85e+004 &     3.5186 \\\hline

160 & 4.44e+004 &     3.7258 & 4.94e+004 &     3.0524 & 4.75e+004 &     3.2595 & 4.68e+004 &     3.6515 \\\hline

180 & 4.29e+004 &     3.8578 & 5.32e+004 &     2.0449 & 5.10e+004 &     2.2668 & 4.57e+004 &     3.7637 \\\hline

200 & 4.18e+004 &     3.9748 & 5.74e+004 &     1.0116 & 5.50e+004 &     1.2419 & 4.49e+004 &     3.8592 \\\hline

220 & 4.09e+004 &     4.0795 & 6.20e+004 &    -0.0045 & 5.93e+004 &     0.2311 & 4.44e+004 &     3.9407 \\\hline

240 & 4.02e+004 &     4.1741 & 6.70e+004 &    -0.9780 & 6.41e+004 &    -0.7394 & 4.40e+004 &     4.0103 \\\hline

260 & 3.96e+004 &     4.2602 & 7.22e+004 &    -1.8951 & 6.90e+004 &    -1.6561 & 4.37e+004 &     4.0695 \\\hline

280 & 3.92e+004 &     4.3388 & 7.77e+004 &    -2.7506 & 7.43e+004 &    -2.5136 & 4.36e+004 &     4.1197 \\\hline

300 & 3.88e+004 &     4.4111 & 8.34e+004 &    -3.5436 & 7.97e+004 &    -3.3102 & 4.35e+004 &     4.1620 \\\hline

500 & {\bf 3.73e+004} &     4.9042 & 1.35e+005 &    -8.5246 & 1.29e+005 &    -8.3127 & 4.47e+004 &     4.2742 \\\hline

\multicolumn{9}{|c|}{$\mu = 5, \tau = 2.5$ } \\\hline

20 & 2.54e+007 &    -2.7911 & 1.10e+023 &  -157.9048 & 8.53e+018 &  -116.7985 & 5.63e+015 &   -84.9895 \\\hline

40 & 4.91e+005 &     3.7130 & 1.37e+040 &  -328.8532 & 8.05e+031 &  -246.5444 & 6.03e+022 &  -155.2917 \\\hline

60 & 4.69e+004 &     4.4065 & 3.59e+057 &  -503.0389 & 1.68e+045 &  -379.7479 & 6.55e+029 &  -225.6465 \\\hline

80 & 3.80e+004 &     4.6991 & 1.29e+075 &  -678.5934 & 4.93e+058 &  -514.4122 & 7.15e+036 &  -296.0278 \\\hline

100 & 3.74e+004 &     4.9027 & 5.53e+092 &  -854.9100 & 1.74e+072 &  -649.8897 & 7.84e+043 &  -366.4253 \\\hline

120 & {\bf 3.73e+004} &     5.0513 & 2.65e+110 & -1031.7135 & 6.92e+085 &  -785.8864 & 8.61e+050 &  -436.8334 \\\hline

140 & 3.73e+004 &     5.1600 & 1.37e+128 & -1208.8552 & 2.99e+099 &  -922.2437 & 9.47e+057 &  -507.2490 \\\hline

160 & 3.74e+004 &     5.2373 & 7.52e+145 & -1386.2455 & 1.37e+113 & -1058.8660 & 1.04e+065 &  -577.6699 \\\hline

180 & 3.76e+004 &     5.2888 & 4.31e+163 & -1563.8262 & 6.62e+126 & -1195.6913 & 1.15e+072 &  -648.0947 \\\hline

200 & 3.78e+004 &     {\bf 5.3182} & 2.55e+181 & -1741.5574 & 3.31e+140 & -1332.6769 & 1.27e+079 &  -718.5224 \\\hline

500 & 4.27e+004 &     4.4523 & Inf       &       -Inf & Inf       &       -Inf & 5.70e+184 & -1775.0426 \\\hline

\end{tabular}
\end{center}}
\end{table}

We also plotted some figures to graphically illustrate the performance of these solvers. Figures (a), (b) and (c) in Figure \ref{fig:compare-deblur} plot the  objective function value versus the iteration number for $\mu=0.1,0.5$ and $1$, respectively. Figures (d), (e) and (f) in Figure \ref{fig:compare-deblur} plot ISNR versus the iteration number for $\mu=0.1,0.5$ and $1$. We did not plot graphs for $\mu=5$, since FaMSA, FaMSA-s and Grad diverged from the very first iteration. From Figure \ref{fig:compare-deblur} we can see the comparisons clearly. Basically, these figures show that FaMSA and FaMSA-s achieve a nearly optimal solution very quickly. We can also see from (b), (c), (e) and (f) that FaMSA-s is always slightly better than FaMSA and MSA is always better than Grad.

We also tested setting $D^{(k)}$ to the identity matrix in MSA and FaMSA, but this choice, as expected, did not give as good results.

\begin{figure}
\centering \subfigure{\includegraphics[scale=0.5]{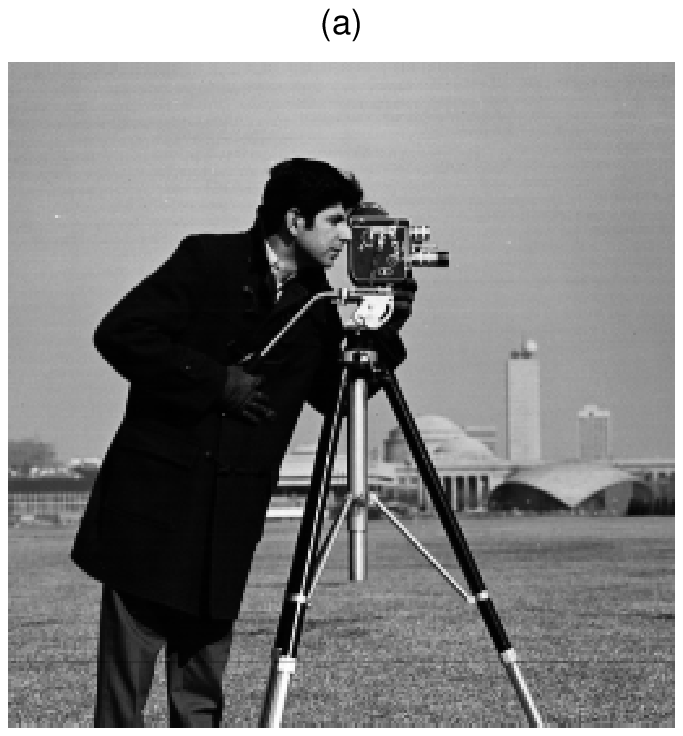}}\hspace{-0.5cm}
\centering \subfigure{\includegraphics[scale=0.5]{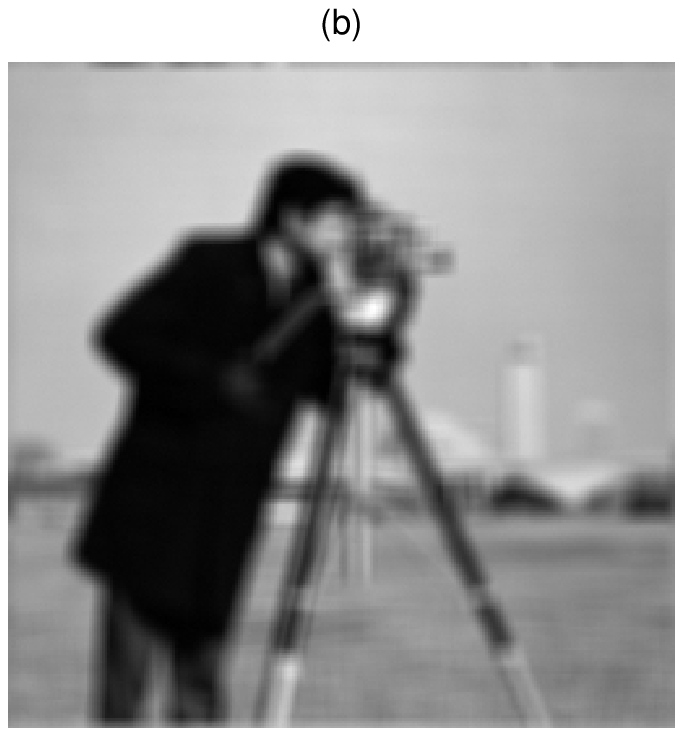}} \hspace{-0.5cm}
\centering \subfigure{\includegraphics[scale=0.5]{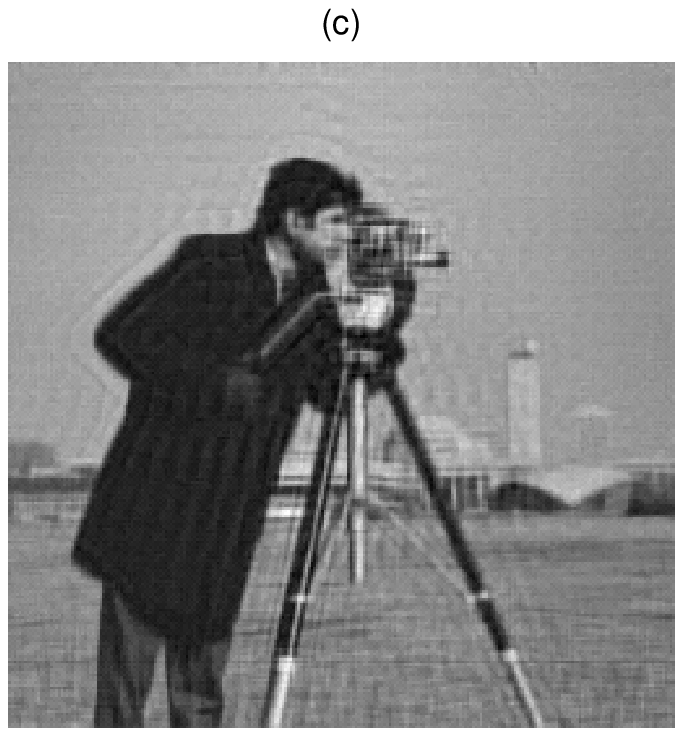}} \hspace{-0.5cm}
\caption{Using MSA to solve \eqref{prob:TV-wavelet-deblur-smooth}. (a): Original image; (b): Blurred image; (c): Reconstructed image by MSA} \label{fig:MSA-deblur}
\end{figure}

To see how MSA performed for the deblurring problem \eqref{prob:TV-wavelet-deblur-smooth}, we show the original (a), blurred (b) and reconstructed (c) cameraman images in
Figure \ref{fig:MSA-deblur}. The reconstructed image (c) is the one that was obtained by applying MSA with $\mu=5$ after 200 iterations. The ISNR of the reconstructed image is 5.3182.
From Figure \ref{fig:MSA-deblur} we see that MSA was able to recover the blurred image very well.

\section{Conclusions}\label{sec:discussion}
In this paper, we proposed two classes of multiple splitting algorithms based on alternating directions and optimal gradient techniques
for minimizing the sum of $K$ convex functions. Complexity bounds on
the number of iterations required to obtain an $\epsilon$-optimal solution for these
algorithms were derived. Our algorithms are all parallelizable, which is attractive for practical applications involving large-scale
optimization problems.

\section*{Acknowledgement}
We would like to thank the anonymous referee for making several very helpful suggestions.

\bibliographystyle{siam}
\bibliography{C:/Mywork/Optimization/work/reports/bibfiles/All}
\end{document}